\documentclass{article}

\usepackage{blindtext} 

\usepackage[sc]{mathpazo} 
\usepackage[T1]{fontenc} 
\usepackage{microtype,amsmath} 

\usepackage[english]{babel} 
 
\usepackage[hmarginratio=1:1,top=32mm,columnsep=20pt]{geometry} 
\usepackage[hang, small,labelfont=bf,up,textfont=it,up]{caption} 
\usepackage{booktabs} 

\usepackage{lettrine} 
\usepackage{subcaption}
\usepackage{enumitem} 
\setlist[itemize]{noitemsep} 
\usepackage{algorithm,algpseudocode}
\usepackage{abstract} 

\usepackage{titlesec} 
\renewcommand\thesection{\Roman{section}} 
\renewcommand\thesubsection{\roman{subsection}} 
\titleformat{\section}[block]{\large\scshape\centering}{\thesection.}{1em}{} 
\titleformat{\subsection}[block]{\large}{\thesubsection.}{1em}{} 

\usepackage{fancyhdr} 
\usepackage{titling} 

\usepackage{hyperref} 
\usepackage{graphicx}

\newtheorem{cor}{Corollary}

\pretitle{\begin{center}\Huge\bfseries} 
\posttitle{\end{center}} 
\title{A velocity tracking approach for the Data Assimilation problem in blood flow simulations} 
\author{%
\textsc{Jorge Tiago, Telma Guerra, Ad\'elia Sequeira}\\
\normalsize CEMAT, Instituto Superior T\'ecnico, Ulisboa \\ 
\normalsize \href{mailto:jftiago@math.tecnico.ulisboa.pt}{jftiago@math.tecnico.ulisboa.pt} 
}
\date{} 
\renewcommand{\maketitlehookd}{%
\begin{abstract}
\noindent 
Several advances have been made in Data Assimilation techniques applied to blood flow modeling. Typically, idealized boundary conditions, only verified in straight parts of the vessel, are assumed. We present a general approach, based on a Dirichlet boundary control problem, that may potentially be used in different parts of the arterial system. The relevance of this method  appears when  computational reconstructions of the 3D domains, prone to be considered sufficiently extended,  are either not possible, or  desirable, due to computational costs. Based on taking a fully unknown velocity profile as the control, the approach uses a discretize then optimize methodology to solve the control problem numerically. The methodology is  applied to a realistic 3D geometry representing a brain aneurysm. The results show that this DA approach may be preferable to a pressure control strategy, and that it can significantly improve the accuracy associated to typical solutions obtained using idealized velocity profiles. 
\end{abstract}
}

\begin{document}

\maketitle


\section{Introduction}\label{intro}

During the last twenty years scientific computing has become an established tool to carry out basic research about the cardiovascular system, in particular its physiopathology. This is due, on the one hand, to the increase of computational power and available medical data and, on the other hand, to the advances in the numerical methods and mathematical models associated to the cardiovascular system. In fact, when put together, these aspects allow us to give a good description of several physiological aspects, as well as some of the pathologies of this complex system. 

Both medical and  scientific communities have also recognized the potential for computational simulations to become an instrument to be used in clinical practice as a complement of diagnosis and even prognosis. Particularly for blood flow, which can be considered as a pulsatile fluid interacting with the vessel wall, mathematical models can be numerically solved. This fact allows to obtain extended measurements not easily acquired in traditional medical devices. 
 Examples of such measurements are the Wall Shear Stress (WSS) or the Oscillating Shear Index (OSI). The WSS and OSI profiles can be indicators of potential  atherosclerosis development or brain aneurysm's rupture, just to name two severe pathologies.  However, diagnosis and prognosis demand  patient-specific and accurate enough simulations, to be considered reliable for medical decisions.

Such a personalization of the results requires the adjustment of both model parameters and boundary conditions, which need to be measured or estimated. {Dealing with this uncertainty remains an active field of research. Simplified models of networks, representing the  cardiovascular system, can be considered reliable to provide information about the flow rate and the pressure average (\cite{Blanco})}. But, whereas these models require relatively few personalized data, detailed spatial distributions of WSS and OSI are more demanding, since the knowledge of the three-dimensional velocity profile describing blood flow becomes determinant. 

A detailed distribution of the WSS or OSI is typically only required locally, in a specific district of an artery.  Therefore, one can consider a surrogate strategy based on the coupling of three levels of detail: a local 3D fluid-structure interaction (FSI) model,  a 1D simplified description for adjacent vessels, and a lumped parameters (0D) model for the remaining vascular tree. This technique, called geometric multiscale (GM), is analyzed in detail by several authors. {We refer to \cite{FNQV}, \cite{Grinberg} and \cite{QVV} for an overview. In these papers, the authors show that GM is a valid approach to obtain a local 3D description, while still capturing the influence of the systemic behavior}. Besides, the procedure can also be fed with patient-specific data,  in order to  make average velocities and pressures more accurate at the level of the simplified networks (\cite{DeVault}, \cite{Ismail}, \cite{Lombardi} and \cite{Pant}). However, when dealing with the artificial boundaries of the domain, where the 3D model is coupled with the surrogate models, the velocity profile is assumed to have an idealized shape, typically constant, or parabolic (see, for instance, \cite{FGNQ} and \cite{KFHJT}). Such assumptions  prevent the velocity profile, inside the region of interest, to capture possible secondary flow and helical effects due to the geometric features  upstream  the district under analysis (\cite{Wk},\cite{Bal}). The lack of precision in the velocity profile -as explained  above- suggests a potential downbeat of the method with regards to WSS or OSI accuracy.

 {Still concerning WSS quantification, in \cite{Osh},  \cite{Castro} and \cite{Jansen}}, significant variability due to changes in the  velocity profile was identified for simulations in brain aneurysms.  As possible workaround, in  \cite{Mz}, \cite{Camp} and \cite{Ml} it was suggested to extend computational domains in order to obtain more accurate profiles. Concerning this aspect, we assume the goal of keeping the full model applied only to the smallest possible realistic domain, mitigating in this way the computational cost associated to  domain extensions. It is worth mentioning that, in a near future, required numerical simulations should couple FSI models with several transport equations describing a clot or plaque formation  inside an artery, which will necessarily be done at the expenses of increasing the  computational cost (\cite{FPS}, \cite{CPM}, \cite{SSST}). Also, image processing still carries several types of uncertainty. To mention but one example, the decision on the smoothing parameters of segmented medical images can lead to an error on the WSS and OSI quantification three times larger than the one resulting from the {uncertainty on model parameters (see \cite{Gambaruto}). Similar conclusion was highlighted in \cite{Cebral}}. This indicates that uncertainty reduction benefits from the choice of smaller domains. 

Once we assume that the computational domain is not big enough to recover the effective real flow structures, additional information must be considered. A natural approach is to consider measurements of the velocity itself and try to adjust the missing boundary conditions. In this way, the measurements can be matched by the computational solution. At this stage, one could ask why not to use the measurements of the three velocity components as the boundary condition itself. In fact, currently it is impossible to obtain velocity measurements covering the entire inlet artificial section, but only at selected points of the vascular domain. We refer to \cite{Morb} and \cite{dReij} for examples of velocity data obtained with phase-contrast MRI. Therefore, regardless of some attempts  made  (see \cite{On} and the works there mentioned), a general strategy, sound  enough to be used in different parts of the cardiovascular system, and which permits to adjust velocity boundary conditions efficiently, remains to be designed.
     
The use of data to improve the computational solution  is the subject of Data Assimilation (DA), which includes different types of approaches and has been used in several fields (see \cite{WZZ} for an overview). The application of DA to hemodynamics has increased in the last decade, mainly for the purpose of model parameter estimation, including material properties needed to properly define FSI models. To name but a few references, we mention \cite{MBXFTCG}, \cite{BV} and \cite{SchClem}. {In \cite{MBXFTCG}  a sequential approach was used to identify stiffness parameters in elastic boundary conditions. Also, in \cite{BV}, a reduced order variational approach was proposed to estimate the Young modulus of the aortic vessel wall. In \cite{SchClem}, a Bayesian analysis was suggested to estimate boundary resistances in the frame of Fontan ventrical palliation surgery.} For an overview on DA approaches in cardiovascular mathematics, we refer to \cite{Mar} or \cite{BDPV}. Several authors also applied DA to the adjustment of boundary conditions. A control approach was suggested in \cite{FVV1} and \cite{FVV2} for the adjustment of the flow rate boundary conditions. In \cite{Heys}  and \cite{RajHeys} the Weighted Least Squares Finite Element Method (WLSFEM) was used to include velocity measurements in the simulations. The method was validated for a recasted form of the Navier-Stokes equations on non-primitive variables. A more flexible technique, based on a variational formulation, was suggested to use velocity measurements in order to adjust pressure boundary values at the artificial boundaries (\cite{Delia3}, \cite{DPV},\cite{DV}). 

In this work, we try to answer some questions concerning the use of the variational formulation. Particularly, we address its possible use as an instrument to adjust 3D velocity profiles on artificial boundaries. Some prospective work considering only idealized velocity profiles, normal to the inlet section, was done in \cite{GTS}. {These idealized profiles were parametrized using up to 2 degrees of freedom. This restrictive assumption resulted in a computational problem with substantially reduced complexity.} Here we abandon this non realistic assumption and we consider - for the first time - the general case, where the velocity profile is allowed to have - for instance- an helical structure. To this end, we use an extended realistic geometry obtained by segmentation of a brain artery with an aneurysm, to generate what we consider the ground truth synthetic data. Afterwards, we truncate the domain into a smaller one, where the ground truth velocity profile is helical. We assume to have velocity data inside several locations in the pathological region, and we apply the variational approach to adjust the boundary velocity profile in order to match the artificially measured data. 

At this stage, we are forced to remain under two non realistic assumptions, which should be dropped in future work. Firstly, the model is assumed to be stationary, in order to neglect the fluid interaction with the vessel walls. Secondly, we assume the velocity data to cover a full section, even though - as already mentioned above- it can only be obtained pointwise. The later assumption allows us to remain within an essentially deterministic frame, and eventually prove that the approach is mathematically sound. If one would have realistic data in mind, a stochastic approach (\cite{DV}, \cite{SchClem}) should be considered. This should be addressed in the future. 

In short, the issues that we try to clarify are the following: can a data assimilation  approach, based on a velocity control problem, be used to obtain a solution that matches measured velocities in a section of the lumen? Will the results improve if more sections are included? Can this approach be preferred to the pressure control problem strategy?

This paper is organized as follows. In Section~\ref{methods} we start by introducing the model for the  blood flow that we will consider; we describe the DA approach and some relevant mathematical issues about it. Then, in Section~\ref{SecNumAp}, we present the numerical algorithm to address our problem, including the  Discretize then Optimize methodology used to solve the resulting control problem. In Section~\ref{results} we present and discuss the numerical results. The results shown include a comparison with the pressure control based approach as well as results in a realistic domain. Finally, in Section~\ref{conclusions}, we end up with several concluding remarks.

\section{Methods}\label{methods} 

The Navier-Stokes equations have been widely accepted as a mathematical model for blood flow in large and medium size arteries (\cite{FQV}). As mentioned before, under the stationary assumption no further interaction with the vessels walls will be considered. Blood flow can also undergo  non-Newtonian characteristics (\cite{RSO1}). Nevertheless, to remain in the frame where mathematical  theory can be directly applied to prove the well posedness of the variational approach, we will consider blood as a Newtonian fluid. The model for the blood flow can read as follows: let the vector function ${\bf u}$ and the scalar function $p$ represent the blood velocity and pressure, respectively.  Both quantities satisfy the momentum and mass balance equations 
\begin{equation}\label{navierstokes}
\left\{ \begin{array}{ll}-\nu\Delta {\bf u}+
	{\bf u}\cdot \nabla {\bf u} +\nabla p= {\bf f}& \qquad
	      \mbox{in }  \Omega,\vspace{2mm} \\
             \nabla \cdot {\bf u}=0&  \qquad \mbox{in }  \Omega,\vspace{2mm}\\
{\bf u}={\bf g}& \qquad \mbox{on }  \Gamma_{in},\\
{\bf u}={0}& \qquad \mbox{on }  \Gamma_{wall},\\
             \nu\partial_n {\bf u}-p{\bf n}={\bf 0}& \qquad \mbox{on }  \Gamma_{out}.
\end{array}\right.\end{equation}

\begin{center}
\begin{figure}[!ht]
\centering
\includegraphics[ scale=0.60]{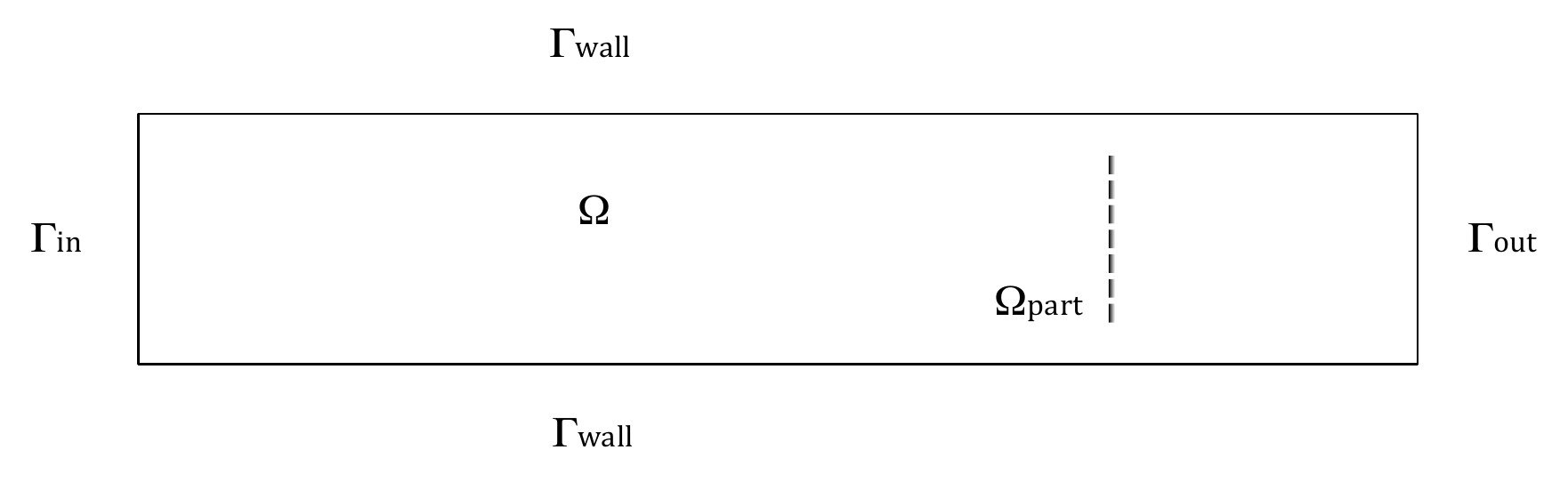}
\centering
\vspace{-0.3cm}
 \caption{\scriptsize{{Example: two dimensional domain.}}}
\label{domain}
 \end{figure}
\end{center}

 Here $\Omega$ represents the vessel domain truncated by  two artificial sections which are set to be the inflow and outflow boundaries, see Figure \ref{domain}. The vector function ${\bf g}$ describes the velocity profile on the inflow boundary $\Gamma_{in}$.  We consider a homogeneous Dirichlet boundary condition on the vessel wall $\Gamma_{wall}$ and a homogeneous Neumann  boundary condition on the outflow boundary $\Gamma_{out}$. The  kinematic viscosity is represented by $\nu$. The body forces are neglected and hence we take ${\bf f} ={\bf 0}$.  

The velocity tracking approach for the DA problem consists of looking for the control function ${\bf g}$ such that  the following cost functional

\begin{equation}\label{costfunctional}
J({\bf u},{\bf g})=\beta_1\int_{\Omega_{part}}|{\bf u}-{ \bf u}_d|^2\,dx+\beta_2\int_{\Gamma_{in}}|\nabla_s {\bf g}|^2\,ds,
\end{equation} 
will be minimized. Here ${\bf u}$ is the solution of (\ref{navierstokes}) corresponding to ${\bf g}$ and  ${\bf u}_d$ represents the data available only on a part of the domain called $\Omega_{part}$. By fixing the parameters $\beta_1$ and $\beta_2$, we decide whether the minimization of $J$ should emphasize a good approximation of the velocity vector to ${\bf u}_d$ or a smoother control measured by the norm of the tangential derivative $\nabla_s(.)$. 

The above problem is a particular case of the broader class of variational problems consisting of different choices for the functional $J$. We remark that in  \cite{Delia3}, \cite{DPV} and \cite{DV} a Neumann control of the type
\begin{equation}\label{neumanncontrol}
[-p{\bf I}+\nu(\nabla{\bf u}+(\nabla{\bf u})^T)]{\bf n}=-g {\bf n}
\end{equation}
was considered at $\Gamma_{in}$.

We will now introduce functional spaces for problem  (\ref{navierstokes}-\ref{costfunctional}).
Let $\Gamma\subset \partial \Omega$ and 
$${\bf H}_{0}^1(\Gamma)= \left\{{\bf v}\in  {\bf L}^2(\Gamma)
        \mid \nabla_s {\bf v} \in  {\bf L}^2(\Gamma), \, \gamma_{\partial\Gamma} {\bf v}={\bf 0} \right\}.$$ 
We constrain the inlet profile as a vector function ${\bf g}\in{\bf \cal U}$ where  
$${\bf \cal U}=\left\{{\bf g}\in  {\bf H}_0^1(\Gamma_{in}):\textrm{such that (\ref{navierstokes}) has a unique weak solution}\right\}. $$ We remark that  ${\bf \cal U}$ is not an empty set as we can take, for instance, ${\bf g}$ such that $\|{\bf g}\|_{H_0^1(\Gamma)}\leq \delta$ for certain $\delta$ small enough (\cite{GST15}).

Now consider $(\Omega_{p_i})_i$ to be a monotone sequence of subsets of $\Omega$, such that 
\begin{equation}\label{subsets}\Omega_{p_1}\subset\Omega_{p_2}...\subset\Omega_{p_m}\subset \Omega.
\end{equation}
In addition, assume also that for all $i\in\{1,...,m\}$, we have $$\partial\Omega_{p_i}=\Gamma_{in}\cup\Gamma_{wall_i}\cup\Gamma_{out_i}$$
where, for all $i\in\{1,...,m\}$, $\Gamma_{out_i}$  are disjoint surfaces corresponding to cross sections of $\Omega${ and } $\Gamma_{wall_i}$ are nonempty wall segments verifying  $\Gamma_{wall_i}\cap\Gamma_{wall}=\Gamma_{wall_i}$.   
Note that the construction of each $\Omega_{p_i}$ in this way  ensures that (\ref{subsets}) is fulfilled, and that  each $\Omega_{p_i}$ itself represents a part of  the vessel $\Omega$. {Therefore, each $\Gamma_{out_i}$ is, in fact, a cross section of $\Omega$.}

We can now state the following {consequence of Theorem 4.5 in \cite{GST15}}:
\begin{cor}\label{existencecontrolsections}

Let $\beta_1, \beta_2>0$ and assume that the data ${\bf u}_d$ is known in a part of the domain  given by $\Omega_{part}=\cup_{i=1}^m{s_i}$ where $s_i=\Gamma_{out_i}$, for all $i\in\{1,...,m\}$. Then there is an optimal solution $({\bf u},{\bf g})\in {\bf H}^1(\Omega)\times {\bf\cal{U}}$ to problem (\ref{navierstokes}-\ref{costfunctional}).
\end{cor}

\subsection{Numerical Approximation}\label{SecNumAp}
In this section we describe the numerical algorithm to solve (\ref{navierstokes}-\ref{costfunctional}). It is based on the  Discretize then  Optimize (DO) approach which  consists of first discretizing the optimal control problem and then solving  the optimization problem (finite dimensional) resulting from the discretization. An alternative approach is the adjoint (indirect) approach, or Optimize then Discretize (OD). {For certain type of parabolic problems, \cite{Betts} and \cite{HinzeTroltzch} indicated that DO approach may be preferred. In \cite{Gunz1} and \cite{Gunz2}, the authors  pointed out that, in nonlinear problems, such us fluid control problems, OD could result in a discrete optimal solution failing to be optimal for the continuous problem.  In \cite{Collis}, in the frame of stabilized advection equations, it was shown that both approaches can lead to different solutions, but, in certain cases, the OD has better asymptotic convergence properties. These conclusions were reinforced in \cite{Heinkenschloss}. Concerning the case of the Navier-Stokes equations, different perspectives were suggested. We refer to \cite{GunzM} for a DO approach in the frame of boundary control, and to \cite{DHinze} for a OD approach in the frame of distributed control. It appears that, at the present stage, no general answer can be given. In particular, concerning problem (\ref{navierstokes}-\ref{costfunctional}), this question remains without unanswered. In \cite{Delia3}, where a pressure type control was considered, the authors obtained better performance of the  DO, in terms of accuracy of the controlled solution. Based on these results,  we have adopted here the DO approach. Nevertheless, a detailed comparison of these two approaches should be the subject of future research.}

Let us assume that we are looking for ${\bf u}\in\mathbf{H}^1(\Omega)$ and for $p\in \mathbf{L}^2(\Omega) $. We consider ${\mathbf{V}=\{\mathbf{v}\in\mathbf{H}^1(\Omega): \ \mathbf{v}|_{\partial \Omega_D}=0\}}$,  where $\partial \Omega_D=\Gamma_{in}\cup \Gamma_{wall}$, $\mathbf{Q}=\mathbf{L}^2(\Omega)$. {Multiplying the first two equations of (\ref{navierstokes}) by test functions $\mathbf{v}\in\mathbf{V}$ and $q\in\mathbf{Q}$, and integrating by parts,  we obtain}
\begin{equation}\label{eq2}
\left\{\begin{array}{ll}
\vspace{0.3cm}
\int\limits_{\Omega}\nu\nabla{\bf u}:\nabla\mathbf{v}+\int\limits_{\Omega}(({\bf u}\cdot\nabla){\bf u})\cdot\mathbf{v}-\int\limits_{\Omega}p\,div\,\mathbf{v}=\int\limits_{\Omega}{\bf f}\cdot\mathbf{v}\\
\int\limits_{\Omega}q\,div\,{\bf u}=0\\
\end{array}
\right.
\end{equation}
{which is the weak form of system (\ref{navierstokes}). The symbol "$:$" represents the inner product of two second-order tensors.}

{To discretize problem (\ref{costfunctional})-(\ref{eq2}) we consider} $\mathbf{V}_h$ and $\mathbf{Q}_h$, subpaces of $\mathbf{V}$ and $\mathbf{Q}$, with finite dimensions $dim(\mathbf{V}_h)=N_u$ and $dim(\mathbf{Q}_h)=N_p$, respectively.  We assume $\mathbf{V}_h$ and $\mathbf{Q}_h$ to represent spaces of Lagrange type Finite Elements, associated to a partition $\tau_h$ of $\Omega$. Therefore the dimensions $N_u$ and $N_p$ tend to infinity when $h$ tends to zero. A map between the nodes describing  $\tau_h$ and the basis functions  with dimensions both $\mathbf{V}_h$ and $\mathbf{Q}_h$ can be defined.

The discretized unknown variables are now given by
\begin{equation}\label{app}
{\bf u}\approx {\bf u}_h=\sum_{j=1}^{N_u}u_j\phi_j\in\mathbf{V}_h, \quad p\approx p_h=\sum_{k=1}^{N_p}p_k\psi_k \in \mathbf{Q}_h
\end{equation}
where $u_j$ and $p_k$ are unknown coefficients to be determined and $\phi_j$ and $\psi_k$ are the shape functions which form a basis of  $\mathbf{V}_h$ and $\mathbf{Q}_h$, respectively. 

Assuming that we can associate some basis functions $(\phi_i)_{i=1...N_o}$ with the nodes in $\Omega_{part}$, and some others to the nodes on  $\Gamma_{in}$, which we refer to as $(\phi_i)_{i=1...N_g}$, we then approximate the  control function in (\ref{costfunctional}) as 
$$\mathbf{g}_h=\sum_{j=1}^{N_g}g_j\phi_j=\sum_{j=1}^{N_g}u_j\phi_j\,.$$
We assume also that ${\bf u}_d$ can be approximated by $${\bf u}_{d,h}=\sum\limits_{i=1}^{N_o}u_{d_i}\phi_i\,.$$  

Let us begin by discretize the cost functional $J$ given by (\ref{costfunctional}). The first term of $J$ becomes 

\begin{eqnarray}\label{1stterm}
&&\int\limits_{\Omega_{part}}\left|\sum_{i}^{N_o}({ u}_i-{{ u}_d}_i)\phi_i\right|^2\,dx
=\int\limits_{\Omega_{part}} \left<\, \sum_{i}^{N_o}({ u}_i-{{ u}_d}_i)\phi_i ,\,\sum_{j}^{N_o}({ u}_j-{{ u}_d}_j)\phi_j \right> \,dx\nonumber\\
&&=\int\limits_{\Omega_{part}} \sum_{i}^{N_o}({ u}_i-{{ u}_d}_i)\sum_{j}^{N_o}({ u}_j-{{ u}_d}_j) \left<\phi_i,\phi_j\right> \,dx\nonumber\\
&&=\sum_{i}^{N_o}({ u}_i-{{ u}_d}_i)\sum_{j}^{N_o}({ u}_j-{{ u}_d}_j)\int\limits_{\Omega_{part}}\phi_i\phi_j\,dx\nonumber\\
&&=({\bf U}-{\bf U}_d)^TM({\bf U}-{\bf U}_d)=\left<({\bf U}-{\bf U}_d),M({\bf U}-{\bf U}_d)\right>\nonumber\\
&&=({\bf U}-{\bf U}_d,{\bf U}-{\bf U}_d)_M=\|{\bf U}-{\bf U}_d \|_{N_o}^2
\end{eqnarray}
where $\|\cdot\|_{N_o}$ is the norm induced by the inner product $(\cdot, \cdot)_{M}$ and ${M}$ is a symmetric $N_o\times N_o$ matrix where each element is given by 
$$m_{ij}=\int\limits_{\Omega_{part}}\phi_i\phi_j\,dx,\,\, i=1,...,N_o,\, j=1,...,N_g.$$ 

For the regularization term we have
\begin{eqnarray}\label{3rdterm}
&&\int\limits_{\Gamma_{in}}\left|\sum_{i}^{N_g}{g}_i\nabla\phi_i\right|^2\,dx=\int\limits_{\Gamma_{in}}\left<\sum_{i}^{N_g}{g}_i\nabla\phi_i, \sum_{j}^{N_g}{g}_j\nabla\phi_j\right>\,dx\nonumber\\
& =&\sum_{i}^{N_g}{g}_i\sum_{j}^{N_g}{g}_j\int\limits_{\Gamma_{in}}\nabla\phi_i:\nabla\phi_j={\bf G}^T{A}{\bf G}\nonumber\\
&=&\left<{\bf G},{A}{\bf G}\right>=({\bf G},{\bf G})_{A}=\|{\bf G}\|^2_{A}
\end{eqnarray}
where $\|\cdot\|_{A}$ is the norm induced by the inner product $(\cdot, \cdot)_{A}$. Matrix ${A}$ is a symmetric $N_g\times N_g$ matrix whose elements are defined by 

$$a_{ij}=\int\limits_{\Gamma_{in}}\nabla\phi_i:\nabla\phi_j\,dx,\,\,i=1,...,N_o,\, j=1,...,N_g,$$ 
where  ":" represents the inner product of two second-order tensors.

{Then, the discretized form of the cost functional (\ref{costfunctional}) becomes:}

\begin{equation}\label{costdis}
J({\bf U}, {\bf G})=\beta_1\|{\bf U}-{\bf U}_d\|^2_{N_o}+\beta_2\|{\bf G}\|^2_A.
\end{equation}

{With respect to system (\ref{eq2}), to deal with convected dominated regimes, a GLS (Galerkin-Least-Squares) stabilization (see \cite{HH}) is adopted here. 
To describe it, let us first consider }
$$a({\bf u}_h,{\bf v}_h)=\int\limits_{\Omega}\nu\nabla{\bf u}_h:\nabla\mathbf{v}_h+\int\limits_{\Omega}(({\bf u}_h\cdot\nabla){\bf u}_h)\cdot\mathbf{v}_h+\int\limits_{\Omega}p_h\,div\,\mathbf{v}_h$$
and 

$$b({\bf u}_h,q_h)=\int\limits_{\Omega}q_h\,div\,{\bf u}_h.$$

{Using this notation, the stabilized version of (\ref{eq2})  consists of finding ${\bf u}_h\in\mathbf{V}_h$ and $p_h\in\mathbf{Q}_h$ such that}

\begin{equation}
\left\{\begin{array}{ll}\label{eq3a}
\vspace{0.3cm}
a({\bf u}_h,{\bf v}_h)+{\cal{L}}^1_h ({\bf u}_h,{\bf f},{\bf v}_h)
=({\bf f},\mathbf{v}_h)\\
b({\bf u}_h,q_h)={\cal{L}}^2_h(p_h,q_h)
\end{array}
\right.
\end{equation}

{where ${\cal{L}}^1_h$ and ${\cal{L}}^2_h$ are defined by}

$${\cal{L}}^1_h ({\bf u}_h,{\bf f},{\bf v}_h)=\sum_{K\in\tau_h}(L({\bf u}_h,p_h)-{\bf f}, \varphi({\bf u}_h,{\bf v}_h))$$
and 

$${\cal{L}}^2_h (p_h,q_h)=(-\frac{1}{\lambda}p_h,q_h)$$ {so that ${\cal{L}}^1_h$ verifies}

\begin{equation}\label{conv}
{\cal{L}}^1_h ({\bf u}_h,{\bf f},{\bf v}_h)=0.
\end{equation}

{Here $\tau_h$ represents a partition of $\Omega$ with characteristic  length $h$, $\lambda$ is a penalty parameter (see \cite{BT}) and $L$ and $\varphi$ are given by }

$$L({\bf u},p)=-\nu\Delta{\bf u}+({\bf u}\cdot\nabla){\bf u}+\nabla p$$

$$\varphi({\bf u}_h,{\bf v}_h)=\delta (({\bf u}_h\cdot\nabla){\bf v}_h+\nu\bigtriangleup{\bf v}_h).$$

{The parameter $\delta$ should be suitably chosen. In this work, the parameter is taken from  \cite{Bazi} (see \cite{Sha} for more details). Notwithstanding, $\delta$ can be optimized in the frame of optimal control problems (\cite{Collis} and \cite{Abr}).}

{In order to obtain the discretization of system (\ref{eq3a}),  let us first discretize the convective term of $a({\bf u}_h,{\bf v}_h)$ and its counterpart in ${\cal{L}}^1_h ({\bf u}_h,{\bf f},{\bf v}_h)$:}

$$\int\limits_{\Omega}(({\bf u}_h\cdot\nabla){\bf u}_h)\cdot\mathbf{v}_h\,+\sum_{K\in\tau_h}\int\limits_{K}(({\bf u}_h\cdot\nabla){\bf u}_h)\cdot\varphi({\bf u}_h,{\bf v}_h).$$

{Using the approximations (\ref{app}) and after some computations, the above expression can be written as}

\begin{align*}
&\left(\sum\limits_{j=1}^{N_u}u_j\sum\limits_{k=1}^{N_u}u_k\int\limits_{\Omega}(\phi_j\cdot\nabla)\phi_k\cdot\phi_i\right)_{i=1,...,N_u}\\
&+\left(\sum_{K\in\tau_h}\sum\limits_{j=1}^{N_u}u_j\sum\limits_{k=1}^{N_u}u_k\int\limits_{K}(\phi_j\cdot\nabla)\phi_k\cdot\left(\delta \left(\sum\limits_{l=1}^{N_u}u_l\phi_l\cdot\nabla{\phi_i}+\nu\Delta{\phi_i}\right)\right)\right)_{i=1,...,N_u}\\
&=({N}({\bf U})+{\cal{N}}({\bf U})){\bf U}
\end{align*}

{where ${\bf U}=(u_1,...,u_{N_u})^T$ and ${{N}}({\bf U})$ and ${\cal{N}}({\bf U})$ are matrices whose elements are defined by}

$$[{{N}}({\bf U})]_{i,j}=\left(\sum\limits_{k=1}^{N_u}u_k\int\limits_{\Omega}(\phi_j\cdot\nabla)\phi_k\cdot\phi_i\right)\quad\forall\, i,j=1,...,N_u$$

$$[{\cal{N}}({\bf U})]_{i,j}=\sum_{K\in\tau_h}\sum\limits_{k=1}^{N_u}u_k\int\limits_{\Omega}(\phi_j\cdot\nabla)\phi_k\cdot\left(\delta \left(\sum\limits_{l=1}^{N_u}u_l\phi_l\cdot\nabla{\phi_i}+\nu\Delta{\phi_i}\right)\right)\quad\forall\, i,j=1,...,N_u.$$

{We now turn our attention to the diffusion term:}

$$\int\limits_{\Omega}\nu\nabla{\bf u}_h:\nabla\mathbf{v}_h+\sum_{K\in\tau_h}\int\limits_{K}-\nu\Delta{\bf u}_h\cdot\varphi({\bf u}_h,{\bf v}_h).$$

{Replacing ${\bf u}_h$ by its corresponding finite approximation we can write}
\begin{align*}
&\left(\nu\sum\limits_{j=1}^{N_u}u_j\int\limits_{\Omega} \nabla\phi_j:\nabla\phi_i\right)_{i=1...N_u}\\
&+\sum\limits_{j=1}^{N_u}{u_j}\sum_{K\in\tau_h}\int\limits_{K}-\nu\Delta\phi_j\cdot\left(\delta\left( \sum\limits_{l=1}^{N_u}u_l\phi_l\cdot\nabla{\phi_i}+\nu\Delta{\phi_i}\right)\right)=
(Q+{\cal Q}){\bf U}, \quad\forall {i=1,...,N_u}
\end{align*}
{where}
$$[Q]_{ij}=\nu\int\limits_{\Omega} \nabla\phi_j:\nabla\phi_i, \quad\forall {i,j=1,...,N_u}
$$
$$[{\cal Q}]_{ij}=\sum_{K\in\tau_h}\int\limits_{K}-\nu\Delta\phi_j\cdot\left(\delta\left( \sum\limits_{l=1}^{N_u}u_l\phi_l\cdot\nabla{\phi_i}+\nu\Delta{\phi_i}\right)\right), \quad\forall {i,j=1,...,N_u}
$$

{Acting similarly for the pressure term, we obtain for the first equation in (\ref{eq3a}),}

$$({Q}+{\cal{Q}}){\bf U}+(N({\bf U)}+{\cal N({\bf U})}){\bf U}+(B^T+{\cal B}){\bf P}={\bf F}$$

{where }

$$[{B}^T]_{i,j}=\int\limits_{\Omega}\psi_j div\,\phi_i\quad\forall\, i=1,...,N_u;\, j=1,...,N_p.$$

$$[{\cal B}]_{ij}=\sum_{K\in\tau_h}\int\limits_{K}\nabla\psi_j\cdot\left(\delta\left(\sum_{l=1}^{Nu}u_l\phi_l\cdot\nabla\phi_i+\nu\Delta\phi_i\right)\right),\quad\forall i=1...N_u; j=1,...N_p.$$

{ as for ${\cal L}_h^2$, we consider}

$$\sum_{E\in\Omega^E}\int\limits_{\Omega^E}\left(\frac{1}{\lambda}p_h,q_h\right)$$

{which, by replacing $p_h$ by it�s corresponding finite approximation, gives}

$$\sum_{k=1}^{N_p}p_k\sum_{K\in\tau_h}\int\limits_{K}\frac{1}{\lambda}\psi_i \psi_k=0\quad\forall i=1,...N_p.$$

{Hence, adding the discretized terms of ${\cal L}_1$ and  ${\cal L}_2$, system (\ref{eq3a}) becomes}
\begin{equation}
\label{dis}
\left\{\begin{array}{ll}
\vspace{0.3cm}
({Q}+{\cal{Q}}){\bf U}+(N({\bf U})+{\cal N}({\bf U})){\bf U}+(B^T+{\cal B}){\bf P}={\bf F}\\
B{\bf U}={\cal B}_1{\bf P} + \textrm{{boundary conditions}}\,\\
\end{array}
\right.
\end{equation}
{where}
$$[{\cal B}_1]_{i,j}=-\sum_{K\in\tau_h}\int\limits_{K}\frac{1}{\lambda}\psi_i \psi_j\quad\forall {i,j=1...N_p}.$$

We remark that vector ${\bf U}=({\bf U}_g,{\bf G})$ includes the controlled velocity coefficients ${\bf G}$ and the uncontrolled ones ${\bf U}_g$ which also depend on ${\bf G}$. Therefore, the stabilized problem can be recast into the general form 

\begin{eqnarray}
\displaystyle\min_{\bf G} F({\bf G})=J({\bf U}({\bf G}),{\bf G}) \label{nlcost} \\ 
C({\bf G})\geq 0 \label{nlconst},
\end{eqnarray}
where (\ref{nlconst}) represents the problem constraints (\ref{dis}),  including boundary conditions.   
In spite of  (\ref{nlcost}-\ref{nlconst}) being finite dimensional, it is a large scale optimization problem with nonlinear constraints and a quadratic cost.  To solve this problem, we use the Sequential Quadratic Programming algorithm, as described in \cite{GMS}. {The algorithm is available in the SNOPT library (\cite{SNO2}) and was tested in several benchmark large scale problems. The iterative procedure requires the evaluation of $F({\bf G})$ which, in turn, implies solving the nonlinear system (\ref{dis}). To solve it, the damped Newton method - as described in \cite{Deu} - was used.} 

We will now briefly describe the algorithm and we refer to \cite{GMS}, for more details.

Let us assume that  the solution ${\bf G}$ of (\ref{nlcost})-(\ref{nlconst}) verifies the Karush-Kuhn-Tucker (KKT) optimality conditions

\begin{align*}
\begin{array}{ll} 
\vspace{0.2cm}
\mathcal{D}C({\bf G})^T\lambda=\mathcal{D}F({\bf G})\\ \vspace{0.2cm}
C({\bf G})^T\lambda = 0\\ \vspace{0.2cm}
C({\bf G}) \geq 0\\
\lambda \geq 0
\end{array}
\end{align*} 
where $\mathcal{D}F$ and $\mathcal{D}C$ are the gradients of $F$ and $C$, respectively, and $\lambda$ is the vector of the Lagrange multipliers. If one is able to find a good initial estimate ${\bf G}_0$ (and corresponding $\lambda_0$), close enough to the optimal ${\bf G}$, the following algorithm produces a sequence that is globally convergent (\cite{GMS}).
\begin{algorithm}
\caption{SNOPT}\label{SNOPT}
  \begin{algorithmic}
 \While{Optimality tolerance of KKT less than threshold} 
    
 \State 1- Determine a quasi-Newton approximation $\mathcal{H}_k$ for the Hessian of the modified Lagrangian
$$\mathcal{L}({\bf G},{\bf G}_k,\lambda_k)=F({\bf G})-\lambda_k^T[C({\bf G}_k)-{\bf G}_k-\mathcal{D}C({\bf G}_k)({\bf G}-{\bf G}_k)].$$

\State 2- Solve the auxiliary Linear Quadratic problem
\begin{align}
&\nonumber \displaystyle\min_C \mathcal{Q}({\bf G},{\bf G}_k,\lambda_k)=F({\bf G}_k)+\mathcal{D}F^T({\bf G}_k)({\bf G}-{\bf G}_k)-\frac{1}{2}({\bf G}-{\bf G}_k)^T \mathcal{H}_k({\bf G}-{\bf G}_k)\\ 
&{\bf G}_k+\mathcal{D}C({\bf G}_k)({\bf G}-{\bf G}_k)\geq 0 \label{lconst}
\end{align}
to obtain the intermediate iterate $(\bar{\bf G}_k,\bar{\lambda}_k,\bar{s}_k)$, where $\hat{s}_k$ is the vector of the slack variables associated to the linear constraints in (\ref{lconst}). 

\State  3- Compute $\alpha_{k+1}\in(0,1]$ as the minimizer of  the merit function
$$M_{\gamma}({\bf G},\lambda,s)=F({\bf G})+\lambda^T(C({\bf G})-s)+\frac{1}{2}\sum_{i=1}^{m}\gamma_i(C_i({\bf G})-s_i)^2$$
 along the line
$$d(\alpha)=({{\bf G}}_k,{\lambda}_k,{s}_k)+\alpha[(\bar{{\bf G}}_k,\bar{\lambda}_k,\bar{s}_k)-({{\bf G}}_k,{\lambda}_k,{s}_k)],$$  
where $s_i$, for $i=1...m$, are the components of $s$ and $\gamma$ is a vector of penalty parameters (see \cite{GMS} for details on how to choose $\gamma$). 

\State 4- Set $({\bf G}_{k+1},{\lambda}_{k+1},{s}_{k+1})=d(\alpha_{k+1})$.

\State 5- Compute the optimality tolerance for the KKT conditions.

\EndWhile

\end{algorithmic}
\end{algorithm}

 \vspace{0.25cm} 
We remark that step 2 of Algorithm \ref{SNOPT}, which concerns the solution of the linear quadratic problem, is conducted using the library SQOPT (\cite{SQOPT}).  

\newpage

\section{Results and Discussion}\label{results}

\subsection{Controlling pressure {\it versus} controlling velocity}
 
As mentioned at the end of Section~\ref{intro}, one of the questions we would like to address is how our approach compares to the approach based on a Neumann control (\cite{Delia3}). To this end, we started by reproducing the results there presented for an idealized 2D straight channel {$\Omega=[0,5]\times[-0.5,0.5]$ with $\Gamma_{in}=\{0\}\times [-0.5,0.5]$ and $\Gamma_{out}=\{5\}\times [-0.5,0.5]$. The observations were assumed to correspond to the sections $\{1\}\times [-0.5,0.5]$, $\{2.5\}\times [-0.5,0.5]$ and $\{4\}\times [-0.5,0.5]$. Taking $\nu=1$, we considered the ground truth solution to be known exactly and given by ${\bf u}=1-4y^2$ (in particular ${\bf u}=0$ on $\Gamma_{wall}=\partial \Omega \setminus (\Gamma_{in}\cup \Gamma_{out})$}. As mentioned above, in \cite{Delia3}, the authors considered the problem of controlling a boundary condition of type (\ref{neumanncontrol}). We refer to this procedure as  solving problem (P2), by opposition to solving problem (\ref{navierstokes}-\ref{costfunctional}), to which we refer as (P1). We remark that while in (P1) the control is a vector  function (2D), in (P2) it is a scalar function. Thus, to solve (P2), the cost function should be properly rewritten, and Algorithm~\ref{SNOPT} may then be  applied in a similar way. 
{In \cite{Delia3}, the weights in the cost function were set to be $\beta_1=\frac{1}{2}$ and $\beta_2=\frac{10^{-9}}{2}$, accordingly to the Morozov Discrepancy Principle associated to a certain fixed signal-to-noise ratio (see, for instance, \cite{ItoK}). Since we were interested in comparing specifically the role of the control nature in the results, we did not include - in this section - any noise on the observations}. 

To solve (P2), as described above, we fixed a Neumann homogeneous condition at $\Gamma_{out}$ and we considered $P_2-P_1$, the usual Lagrange linear FEM  corresponding to 27K degrees of freedom {(maximum element size $h=1/20$)} for the velocity. The assembly of the FEM matrices required to obtain  the equivalent to system (\ref{dis})  was done with COMSOL Multiphysics (\cite{C1}). {Since, at this stage, the Reynolds number was very small, the matrices corresponding to the convective and stabilization terms were neglected}. To solve the linearized systems at the iteration level, the PARDISO library was used. 
{The result gave a controlled solution that approximates the exact solution with a relative error of 0$.00112$, that is, of order $\approx 0.1 \%$. Correspondingly to the conclusions in \cite{Delia3}, this means that the Neumann control was able to successfully adjust the solution to the data. For this reason, we used this percentage as the reference relative error to fix the weights $\beta_1$ and $\beta_2$ in our comparative example - which will be described next.} 

\begin{figure}
\centering
\begin{minipage}{0.45\textwidth}
\includegraphics[scale=0.4]{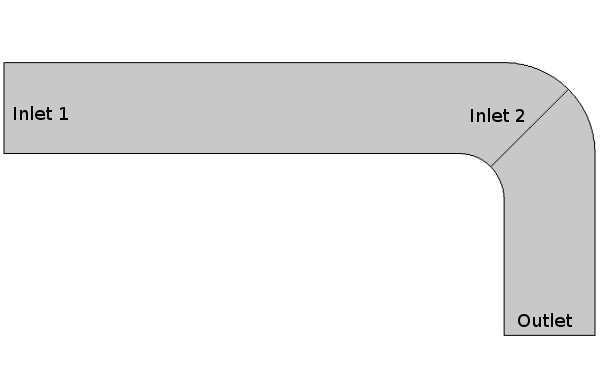}\hspace{.5cm}
\end{minipage}
\begin{minipage}{0.47\textwidth}
\includegraphics[scale=0.35]{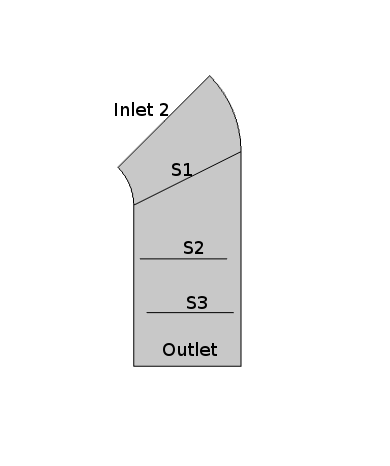}
\end{minipage}
\caption{Ground truth domain (left); Working domain $\Omega$ and $\Omega_{part}=S_1\cup S_2\cup S_3$ (right).}
\label{2dchannel}
\end{figure}

We considered the previous domain and extended it to obtain the  {curved vessel represented in Figure~\ref{2dchannel} (left)}. We will refer to this extended domain as  the ground truth domain. {As it is well known (\cite{QVV}), in a straight channel, even for at physiological Reynolds numbers, pressure contours tend to remain parallel to the cross sections of the computational domain. In a curved vessel, however, that is no longer true, even for Reynolds numbers bellow typical physiologic values.} 

\begin{figure}
\centering
\begin{minipage}{0.45\textwidth}
\includegraphics[scale=0.4]{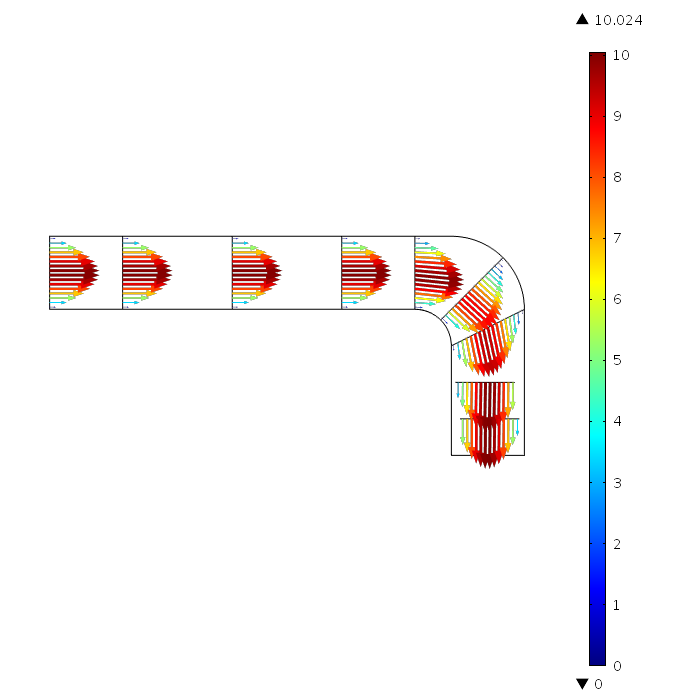}\hspace{0.5cm}
\end{minipage}
\begin{minipage}{0.45\textwidth}
\includegraphics[scale=0.4]{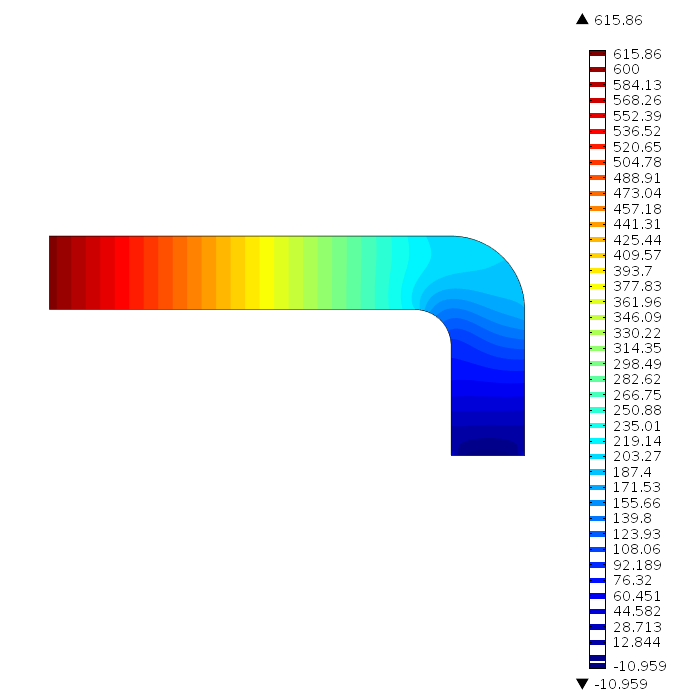}
\end{minipage}
\caption{Ground truth  solution. On the left: velocity magnitude (m/s). On the right: pressure (Pa).}
\label{Ground2dresults}
\end{figure}

{To illustrate this, fixed the same model as before, but we slightly increased the Reynolds number by considering a parabolic profile $10(1-4^2)$ at {\it inlet 1}. We then solved system (\ref{navierstokes}) to obtain the ground truth solution ${\bf u}_d$. 
As the Reynolds number was now higher, we considered all the terms in system (\ref{dis}), including the stabilizing terms. An unstructured  mesh corresponding to 43K degrees of freedom (max $h=1/20$) was used}.  The nonlinear system was solved using the damped  Newton's method, as mentioned in Section~\ref{methods}. {The ground truth solution is represented in Figure~\ref{Ground2dresults}}. We can see that the pressure contours are no long parallel to cross sections within the curve, and the velocity profile loses the parabolic shape on those cross sections.

\begin{figure}
\centering
\begin{minipage}{0.4\textwidth}
\includegraphics[scale=0.32]{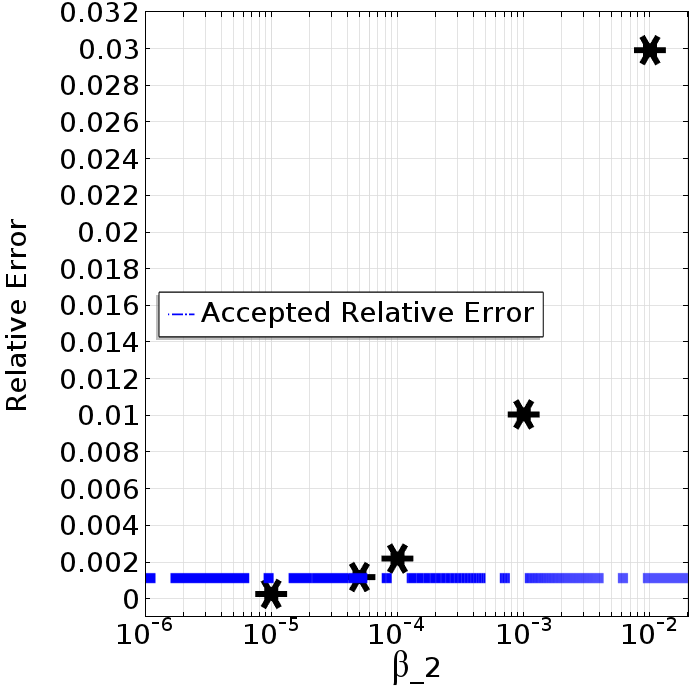}
\end{minipage}
\begin{minipage}{0.45\textwidth}
\includegraphics[scale=0.32]{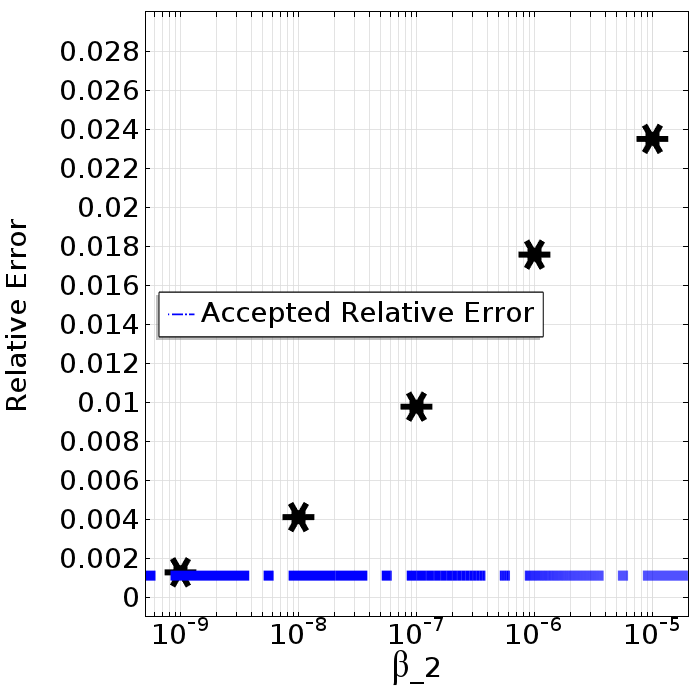}
\end{minipage}
\caption{Illustration of results in Table~\ref{tableRelErrors}: Relative errors $RE_{\Omega_{part},\beta_2}$ for (P1) (left) and (P2) (right).}
\label{PlotRelErrors}
\end{figure}

{Our aim next was to mimic a more realistic situation, where the unknown inlet boundary condition did not correspond to a parabolic velocity profile, normal at $\Gamma_{in}$, nor to a pressure profile, that could be assumed axial dependent.} For this reason, we truncated the channel at the section labeled {\it inlet~2}, which became the new artificial inlet  of the shorter domain represented in Figure~\ref{2dchannel}(right). We call this domain $\Omega$ and its inlet boundary $\Gamma_{in}$. Therefore, we put ourselves into the scenario where we would like to fix a boundary condition at $\Gamma_{in}$ so that the solution in the shorter domain $\Omega$ would match, as much as possible, the ground truth solution ${\bf u}_d$. For the observations,  we assumed to have measured exactly the  velocity profiles of the true solution ${\bf u}_d$ at $\Omega_{part}=S1\cup S2 \cup S3$, where $S1$, $S2$ and $S3$ are lines that were chosen arbitrarily inside $\Omega$ (Figure~\ref{2dchannel}, left). {Before solving both problems (P1) and (P2), we needed to set $\beta_1$ and $\beta_2$. Following what was done above, concerning the example in \cite{Delia3}, we fixed $\beta_1=\frac{1}{2}$ and we looked for $\beta_2$ so that a relative error 
$$RE_{\Omega_{part},\beta_2}=\frac{\|{\bf u}_{\beta_2}-{\bf u}_d\|_{\mathbf{L}^2(\Omega_{part})}}{\|{\bf u}_d\|_{\mathbf{L}^2(\Omega_{part})}},$$ 
verifies $RE_{\Omega_{part},\beta_{2}}\approx 0.00112$. In the expression of the relative error, ${\bf u}_{\beta_2}$ represents the solution of the control problem associated to $\beta_2$. We did this by heuristically fixing a sample for $\beta_2$ and evaluating the corresponding relative errors obtained using  Algorithm~\ref{SNOPT} with an optimality tolerance of $10^{-6}$. The results are shown in  Figure~\ref{PlotRelErrors} and Table~\ref{tableRelErrors}.}

\begin{table}[ht!]
\caption{Relative errors $RE_{\Omega_{part},\beta_2}$, for both (P1) and (P2) approaches. }
\centering
%
\begin{tabular}{cccc}
\toprule
  $\beta_2$  & $(P1)$\\
\midrule
$0.5\times{10^{-2}}$ &  $0.02992$ \\
$0.5\times{10^{-3}}$ &  $0.01005$ \\
$0.5\times{10^{-4}}$&  $0.00219$ \\
$2.5\times{10^{-5}}$&  $0.00118$ \\
$0.5\times{10^{-5}}$&  $0.000259$ \\
\bottomrule
\end{tabular}
\begin{tabular}{cccc}
\toprule
$\beta_2$  & $(P2)$\\
\midrule
$0.5\times{10^{-5}}$ & $0.02350$\\
$0.5\times{10^{-6}}$ & $0.01757$\\
$0.5\times{10^{-7}}$ & $0.00978$\\
$0.5\times{10^{-8}}$ & $0.00412$\\
$0.5\times{10^{-9}}$ & $0.00129$\\
\bottomrule
\end{tabular}

\label{tableRelErrors}
\end{table}

{From these conclusions we fixed $\beta_2=2.5\times{10^{-5}}$ for (P1) $\beta_2=0.5\times{10^{-9}}$ for (P2).} 

\begin{figure}
\centering
\begin{subfigure}[b]{0.4\textwidth}
\includegraphics[scale=0.24]{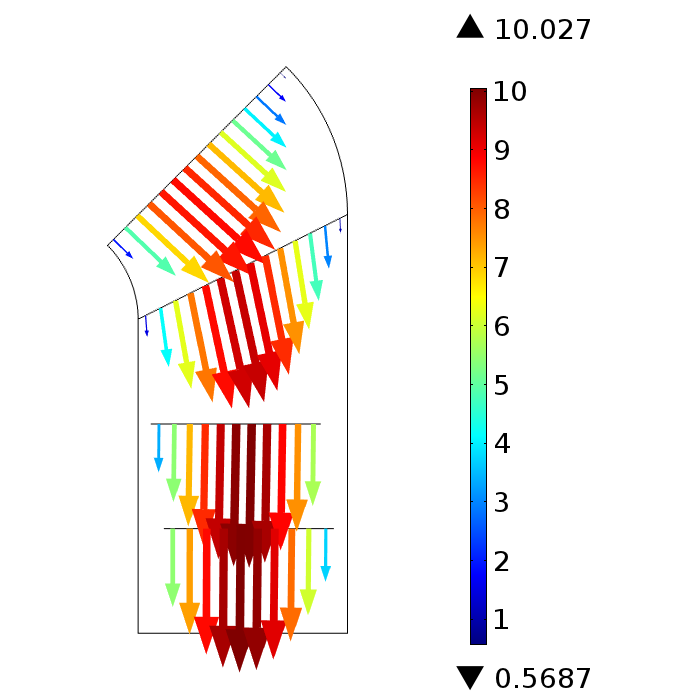}
\caption{Ground truth velocity}
\end{subfigure}    
\begin{subfigure}[b]{0.4\textwidth}
\includegraphics[scale=0.24]{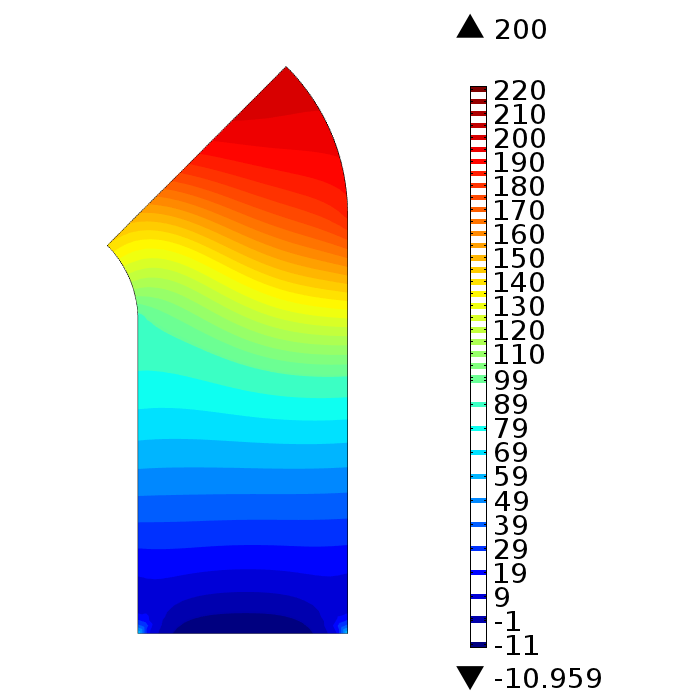}
\caption{Ground truth pressure}
\end{subfigure}

\begin{subfigure}[b]{0.4\textwidth}
\includegraphics[scale=0.24]{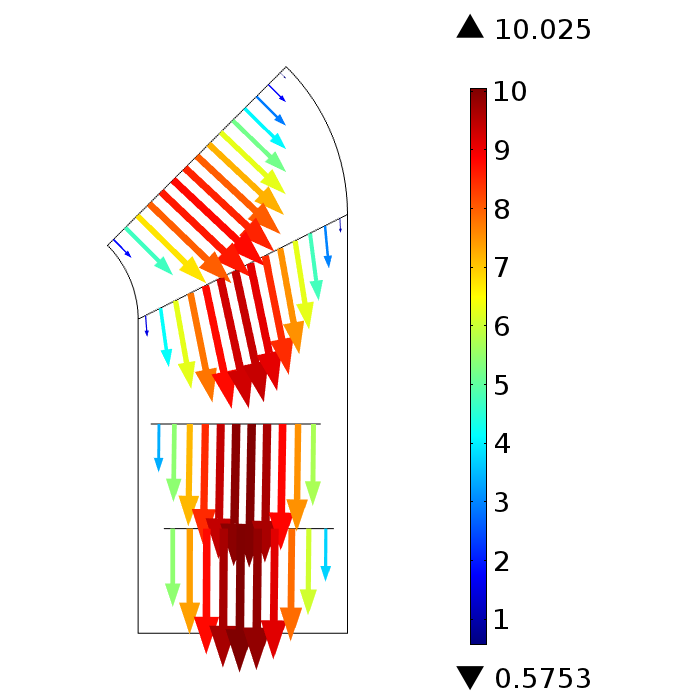}
\caption{Velocity for (P1)} 
\end{subfigure}
\begin{subfigure}[b]{0.4\textwidth}
\includegraphics[scale=0.24]{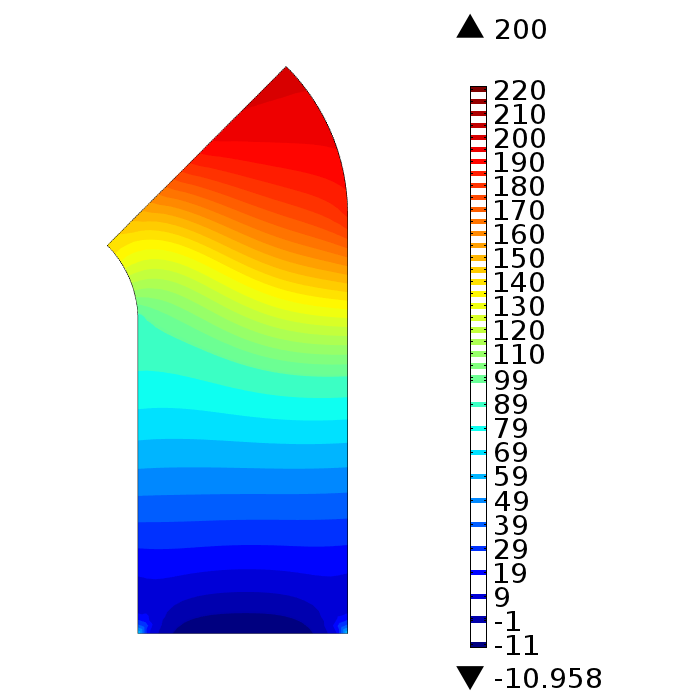}
\caption{Pressure for (P1)}
\end{subfigure}

\caption{ First row: ground truth velocity magnitude (m/s) (left) and ground truth pressure (Pa) (right).  Second row: controlled solution  (P1) - velocity magnitude (m/s) (left) and  pressure (Pa) (right).}
\label{2dresultsvelcont}
\end{figure}

\begin{figure}
\centering
\begin{subfigure}[b]{0.4\textwidth}
\includegraphics[scale=0.24]{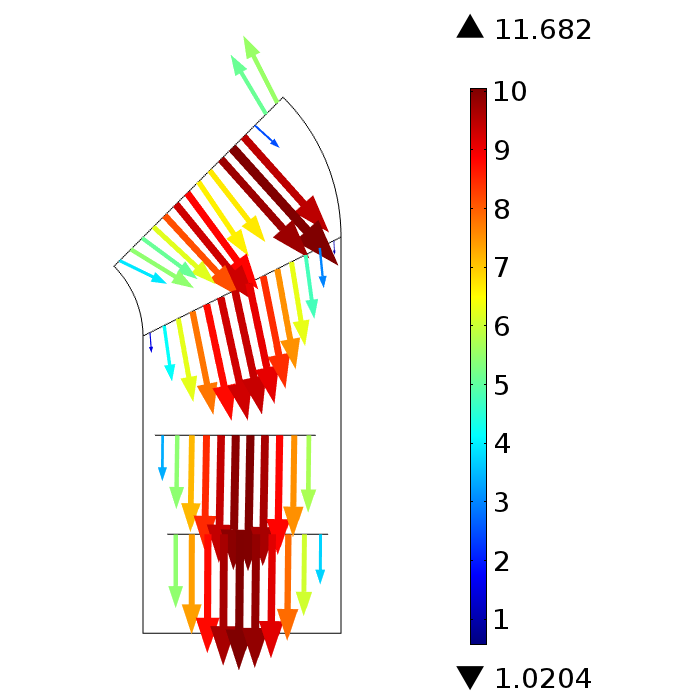}
\caption{.}
\end{subfigure}
\begin{subfigure}[b]{0.4\textwidth}
\includegraphics[scale=0.24]{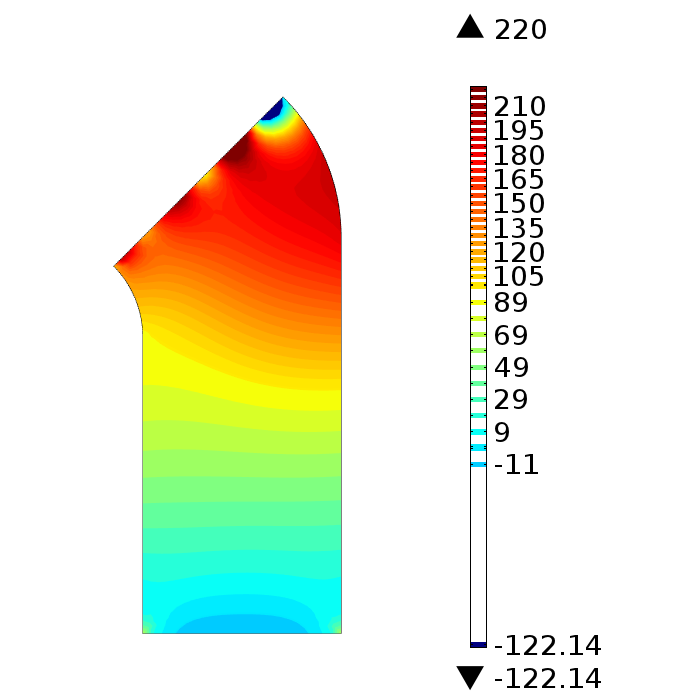}
\caption{.}
\end{subfigure}

\begin{subfigure}[b]{0.4\textwidth}
\includegraphics[scale=0.24]{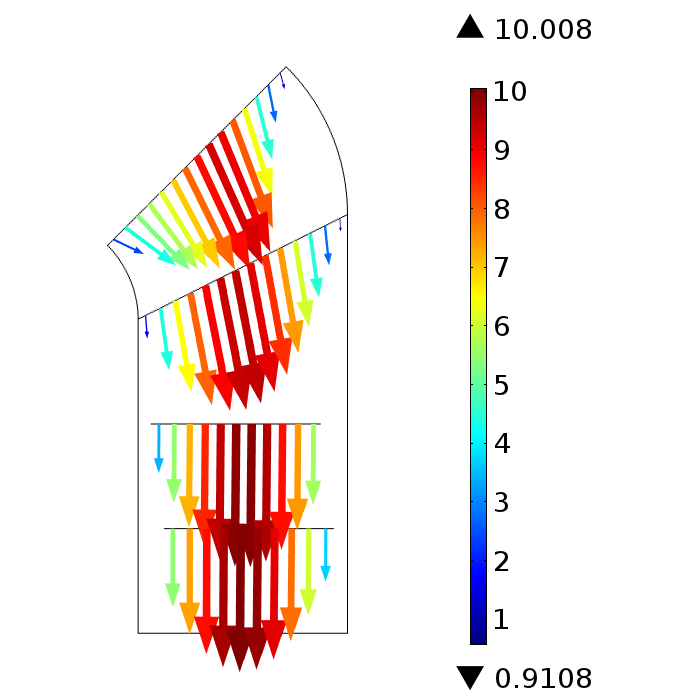}
\caption{.}
\end{subfigure}
\begin{subfigure}[b]{0.4\textwidth}
\includegraphics[scale=0.24]{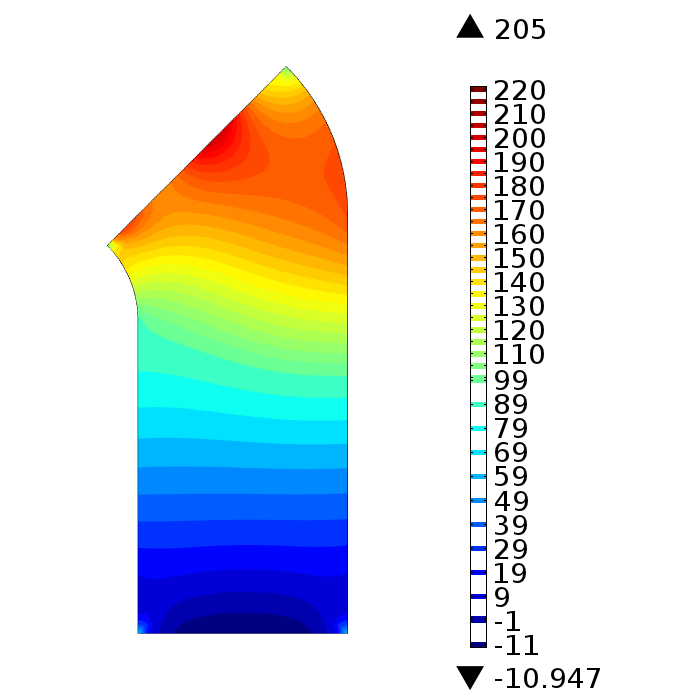}
\caption{.}
\end{subfigure}
\caption{ First row:  controlled solution obtained with (P2) - velocity magnitude (m/s) (left) and  pressure (Pa) (right). Bottom row: same results obtained with $\beta_2=0.5\times 10^{-6}$ (P2r).}
\label{2dresultsprescont}
\end{figure}

{The solutions obtained for (P1) and (P2) are represented in the second row of Figure~\ref{2dresultsvelcont}, and first row of Figure~\ref{2dresultsprescont}, respectively.}

The results show that the solution obtained with the velocity control is qualitatively closer to the ground truth solution, represented in the first row of Figure~\ref{2dresultsvelcont}. To quantify these different performances, we use the relative error of the controlled solutions with respect to ${\bf u}_d$, evaluated at different sites. 
In Table \ref{table 2d} we present the values for  
 $$RE_{\Omega}=\frac{\|{\bf u}-{\bf u}_d\|_{\mathbf{L}^2(\Omega)}}{\|{\bf u}_d\|_{\mathbf{L}^2(\Omega)}},$$
where $\|\cdot\|_{\mathbf{L}^2(\Omega)}$ is the $\mathbf{L}^2(\Omega)$ norm, and for $RE_{\Gamma_{in}}$ and $RE_{\Omega_{part}}$, which are computed analogously. We also indicate the final value for the cost functional and the number of cost  evaluations. It can be seen that, while a relative error on the observations site is kept on the same order, the solution of (P1) is globally more near to ${\bf u}_d$ than the solution of (P2). {Actually, looking closer to the later pressure profile (Figure~\ref{2dresultsprescont}, 1st row), some oscillations can be seen at the inlet. This indicates that, although the relative error on the observations was of order $0.1\%$, the weight $\beta_2=0.5\times 10^{-9}$ almost neglected the regularizing effect of the second term in the cost function. An increase in $\beta_2$ improves the regularizing effect, but at the expenses of distancing from the desired relative error. We illustrate this by considering the case (P2r) with $\beta_2=0.5\times 10^{-6}$, for which the results are shown in  Table~\ref{table 2d} and in the second row of Figure~\ref{2dresultsprescont}.}

\begin{table}[ht!]
\caption{Relative errors, final value for cost functional (J), number of objective evaluations NE for both (P1) and (P2) approaches. }
\centering

\begin{tabular}{cccccc}
\toprule
Approach  & $RE_{\Omega}$ &$RE_{\Gamma_{in}}$ & $RE_{\Omega_{part}}$ & Cost&NE\\
\midrule
(P1)    & $0.00517$ & $0.02286$ & $0.00118$ & $0.00953$&$230$\\
(P2)   & $0.13443$ & $0.57843$ & $0.00129$ & $7.40043e-4$&$126$\\
(P2r)   & $0.10068$ & $0.47189$ & $0.01757$ & $0.03435$&$68$\\
\bottomrule
\end{tabular}
\label{table 2d}
\end{table}

 These results indicate that, if the accuracy of the numerical solution is our aim,  in some scenarios - like the one just illustrated - a velocity control approach can be more convenient. A higher number of cost evaluations might however be required. Such conclusion does not invalidates the fact that (P2) is prone to perform well, when pressure contours align with cross section in the region close to the inlet, as shown in \cite{Delia3} and \cite{DPV}.

\begin{figure}
\centering
\begin{minipage}{0.4\textwidth}
\includegraphics[scale=0.32]{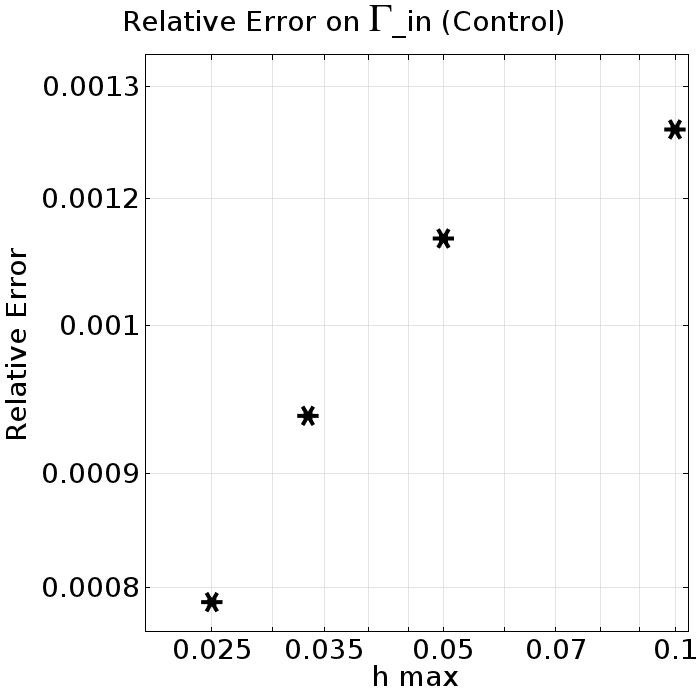}
\end{minipage}
\begin{minipage}{0.45\textwidth}
\includegraphics[scale=0.32]{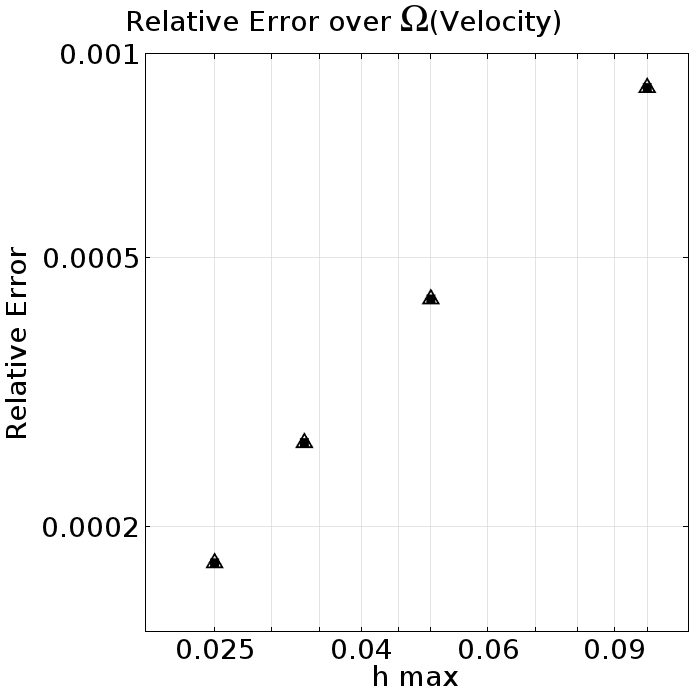}
\end{minipage}
\caption{Mesh convergence: Relative errors for the control vector, on $\Gamma_{in}$ (left) and velocity vector on the controlled domain (right). $h_{max}=1/40$, $h_{max}=1/30$, $h_{max}=1/20$ and $h_{max}=1/10$. }
\label{PlotMeshConvRelErrors}
\end{figure}

{Next, we assessed the convergence of the solution pair $(\mathbf{u},\mathbf{g})$ with respect to mesh refinement. We considered the desired solution as the one obtained by solving (P1) with a maximum mesh size of $h=1/80$. Then, we computed the relative errors both control and velocity solutions, which we represent in Figure~\ref{PlotMeshConvRelErrors}. The convergence test resulted positive.}

{Until now, (P1) was only solved for small Reynolds numbers. This allowed us to illustrate the advantages of velocity control, while remaining close to the example analyzed in \cite{Delia3}. These Reynolds numbers were still far for common physiologic values found in the cardiovascular system. To explore the robustness with respect to such values, we solved (P1) for increasing Reynolds numbers. The results are presented in Table~\ref{table Reynolds}. We can see that the  number of iterations, required to solve Algorithm~\ref{SNOPT}, increases significantly with the Reynolds number. Also, since we kept the weights $\beta_1$ and $\beta_2$ fixed, there is an increase of the relative error on the observations site $\Omega_{part}$. Naturally, this error propagates to the rest of the solution. To keep the relative error within acceptable values, let say, around $0.1\%$, we need to adjust the weight parameters, as it was done above. To illustrate this, we show the results obtained by considering $\beta_2=0.1\times 10^{-5}$ in Table~\ref{table ReynoldsII}. We can see that with these parameters the relative errors can be of order $0.1\%$.}

\begin{table}[ht!]
\caption{Relative errors, number of SNOPT iterations for $\beta_2=2.5\times{10^{-5}}$.}
\centering

\begin{tabular}{cccccc}
\toprule
Reynolds  & $RE_{\Omega}$ &$RE_{\Gamma_{in}}$ & $RE_{\Omega_{part}}$ & &Iterations\\
\midrule
$6.67$   & $0.00517$ & $0.02286$ & $0.00118$ & &$230$\\
$100$    & $0.00836$ & $0.04528$ & $0.00160$ & $$&$260$\\
$200$   & $0.01266$ & $0.07051$ & $0.00214$ & $$&$487$\\
$300$   & $0.01465$ & $0.08328$ & $0.00246$ & $$&$498$\\
$400$   & $0.01587$ & $0.09151$ & $0.00267$ & $$&$454$\\
\bottomrule
\end{tabular}
\label{table Reynolds}
\end{table}
\begin{table}[ht!]
\caption{Relative errors, number of SNOPT iterations for $\beta_2=0.1\times{10^{-5}}$. }
\centering

\begin{tabular}{cccccc}
\toprule
Reynolds  & $RE_{\Omega}$ &$RE_{\Gamma_{in}}$ & $RE_{\Omega_{part}}$ & &Iterations\\
\midrule
$100$    & $0.00634$ & $0.0372$ & $0.000812$ & $$&$374$\\
$200$   & $0.00983$ & $0.05905$ & $0.00117$ & $$&$680$\\
$300$   & $0.01141$ & $0.07013$ & $0.00138$ & $$&$760$\\
$400$   & $0.01243$ & $0.0775$ & $0.0015$ & $$&$774$\\
\bottomrule
\end{tabular}
\label{table ReynoldsII}
\end{table}

%
%

\subsection{The velocity control DA approach applied to a realistic domain}\label{realistic}
In this  section we present the numerical results found after applying  the DA  approach (P1) to a realistic geometry obtained from  the segmentation of Computed Tomography (CT) data sets of a saccular brain aneurysm.

\begin{figure}
\centering
\begin{minipage}{0.45\textwidth}
\includegraphics[scale=0.45]{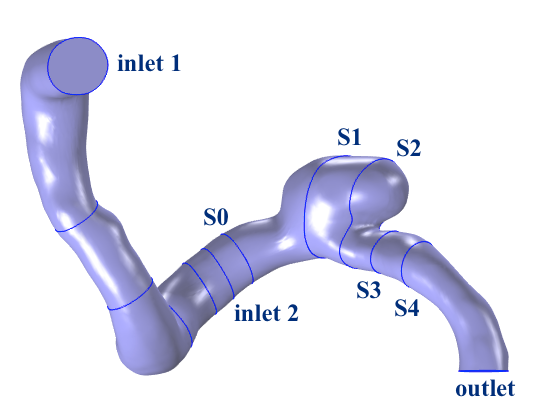}
\end{minipage}
\begin{minipage}{0.45\textwidth}
\includegraphics[scale=0.45]{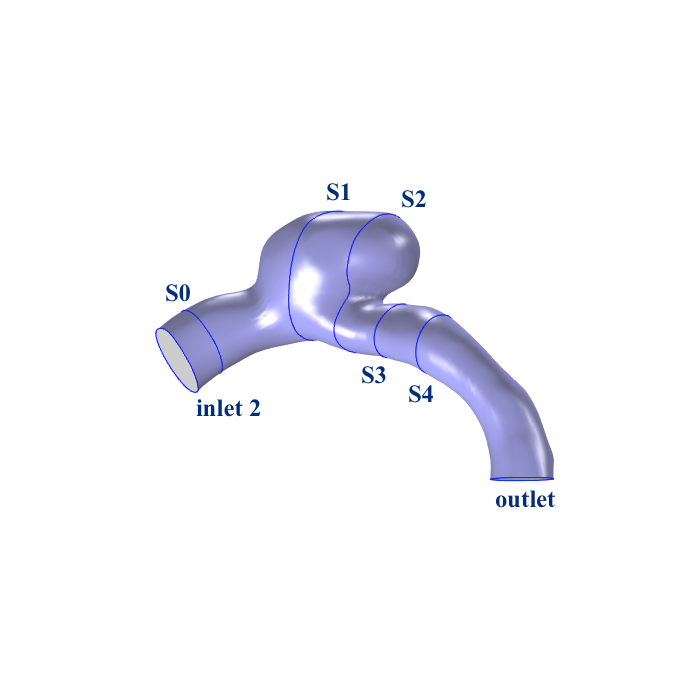}
\end{minipage}
\caption{Ground truth geometry (left); Working geometry $\Omega$ (right).}
\label{aneumesh}  
\end{figure}

As in the previous example, we used  an extended ground truth domain to computed the ground solution ${\bf u}_d$. This solution was used both to select the measured data and to estimate the accuracy of the method. The ground truth domain is represented in Figure~\ref{aneumesh} (left). For the model parameters, we considered $\nu=\frac{\mu}{\rho}$ with $\mu=3.67\times 10^{-3}\, Pa.s$, a value within the range suggested in \cite{Gambaruto}. Also, we took  $\rho=1050\,Kg/m^3$ and we fixed a laminar inflow profile - normal to the inlet -  which corresponds to a flow rate of $Q=4\times 10^{-6}m^3/s$. Again, these are typical parameters used for blood flow simulations (see \cite{Gambaruto}). At the inlet this values implied a physiological Reynolds number of $367$. No slip boundary conditions were imposed on the vessel wall and zero normal stress (Neumann homogeneous) was fixed on the outflow boundary. To obtain system~\ref{dis} and for its numerical solution we adopted the same choices as in the previous example, except that P1-P1 finite elements were used instead. We remark that the system was stabilized using GLS.

{ First, we analyzed the case when the same degrees of freedom were used both to generate ${\bf u}_d$ and for  the DA procedure. Subsequently, the ground truth solution was generated using a finer mesh. The second scenario was considered in order to avoid the so called {\it inverse crime} problem. This is the case when the same model and discretization are used both to generate the synthetic data - from where the observations are chosen - and to solve the control (inverse) problem. Indeed, some particular inverse crime problems can have a trivial solution (\cite{KS}). Although this is not necessarily the case for fluid control problems (see, for instance, \cite{Gunz1} and \cite{Gunz2} for some simple examples with non-trivial solution), we still distinguish between both scenarios. Thus, we can change the observed data either by generating it using a different (finer) mesh, adding noise to it, or both, as we have done it this study.} 

\subsubsection{DA problem - the ``inverse crime'' situation.}

\begin{figure}
\centering   
\includegraphics[scale=0.4]{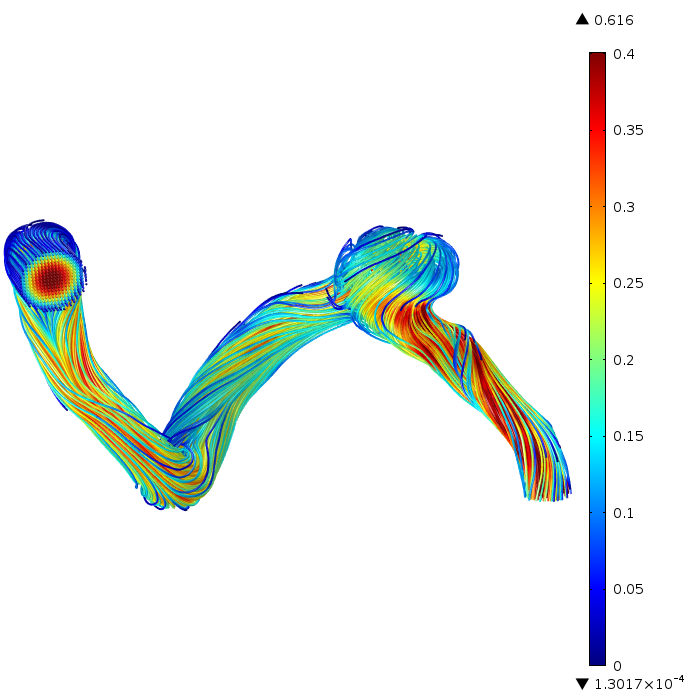}	
\caption{Streamlines representation of the ground truth solution ${\bf u}_d$.}
\label{datastream}
\end{figure}

The first result concerns the case where ${\bf u}_d$ is generated using 213K degrees of freedom. The streamlines of the ground truth solution can be seen in Figure~\ref{datastream}. The helical structures downstream the first steep curvature are evident.

We consider $\Omega$ to be the subdomain starting in section {\it inlet  2}, which is identified in Figure~\ref{aneumesh}, on the right. We identify this section with $\Gamma_{in}$ in problem  (\ref{navierstokes})-(\ref{costfunctional}). We set the goal of  finding a velocity boundary condition to use at this section in such a way that the corresponding solution matches ${\bf u}_d$. Additionally, we assume to have exact measurements of the velocity on $\Omega_{part}=S1\cup S2\cup S3\cup S4$ where $S1$ $S2$, $S3$ and $S4$ are the  sections represented in Figure~\ref{aneumesh}. The later assumption, concerning the exactness of the measurements, will be relaxed in the next study case. 

Concerning the choice of the weights for the cost function here we fixed $(\beta_1,\beta_2)=(10^{5},10^{-3})$. This choice will be justified in section~\ref{secnocrime}, where the presence of noise on the observations will be considered.

\begin{figure}
\centering
 \begin{subfigure}[b]{0.3\textwidth}
\includegraphics[scale=0.35]{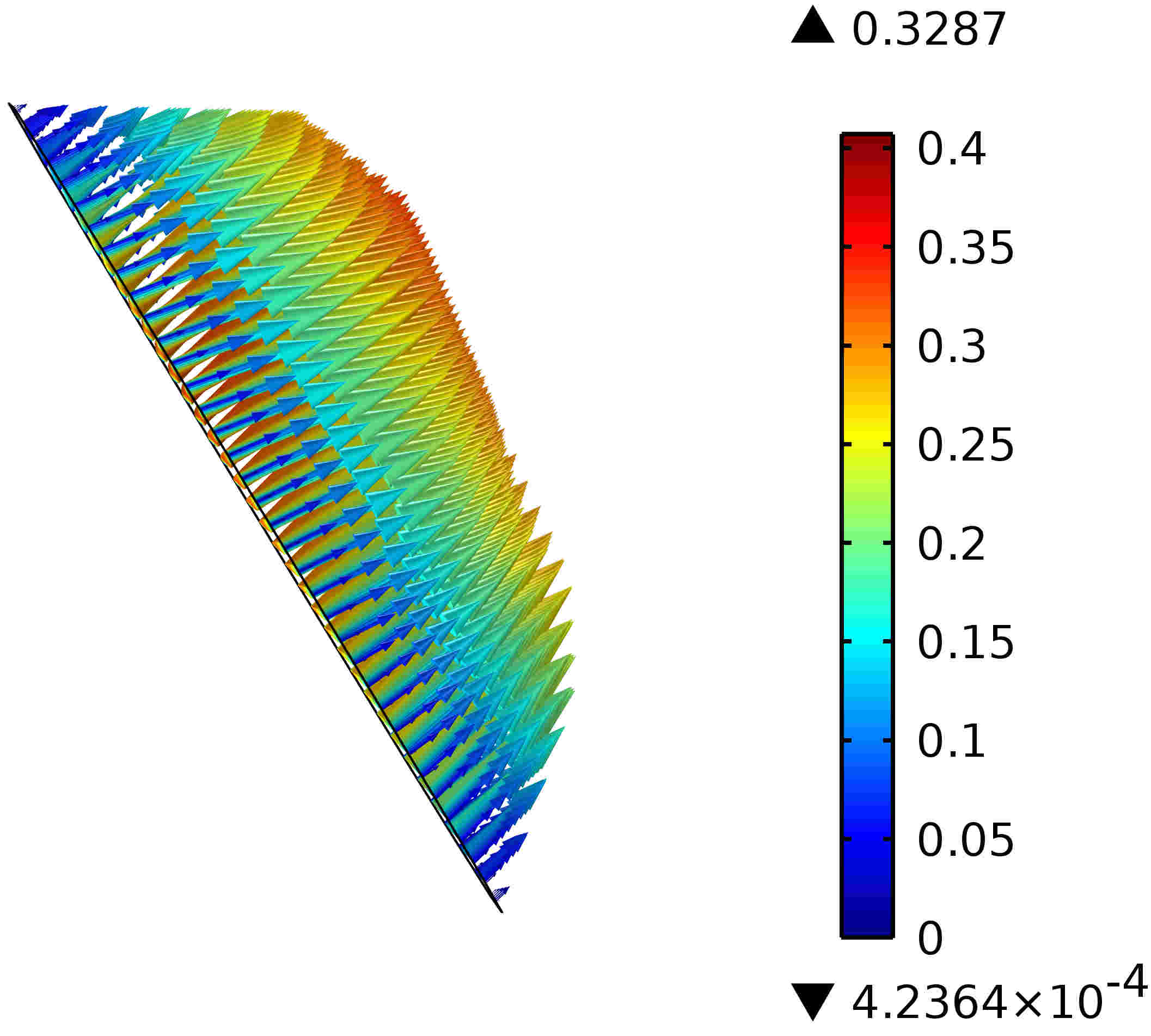}
\caption{Ground truth   ${\bf u}_d$.} 
\end{subfigure}
 \begin{subfigure}[b]{0.3\textwidth}
\includegraphics[scale=0.35]{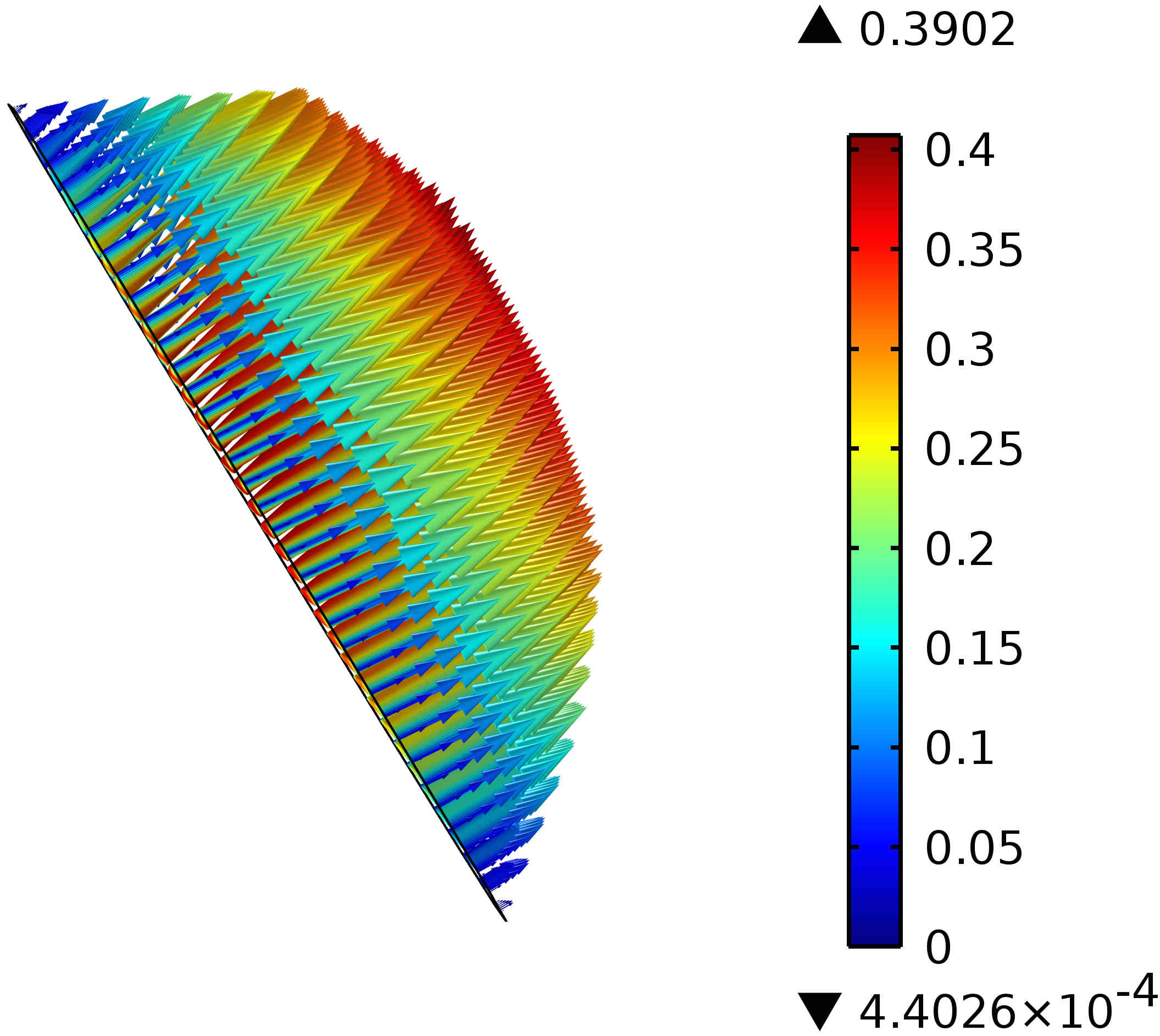}
\caption{Controlled   ${\bf u}$.}
\end{subfigure}
 \begin{subfigure}[b]{0.3\textwidth}
\includegraphics[scale=0.35]{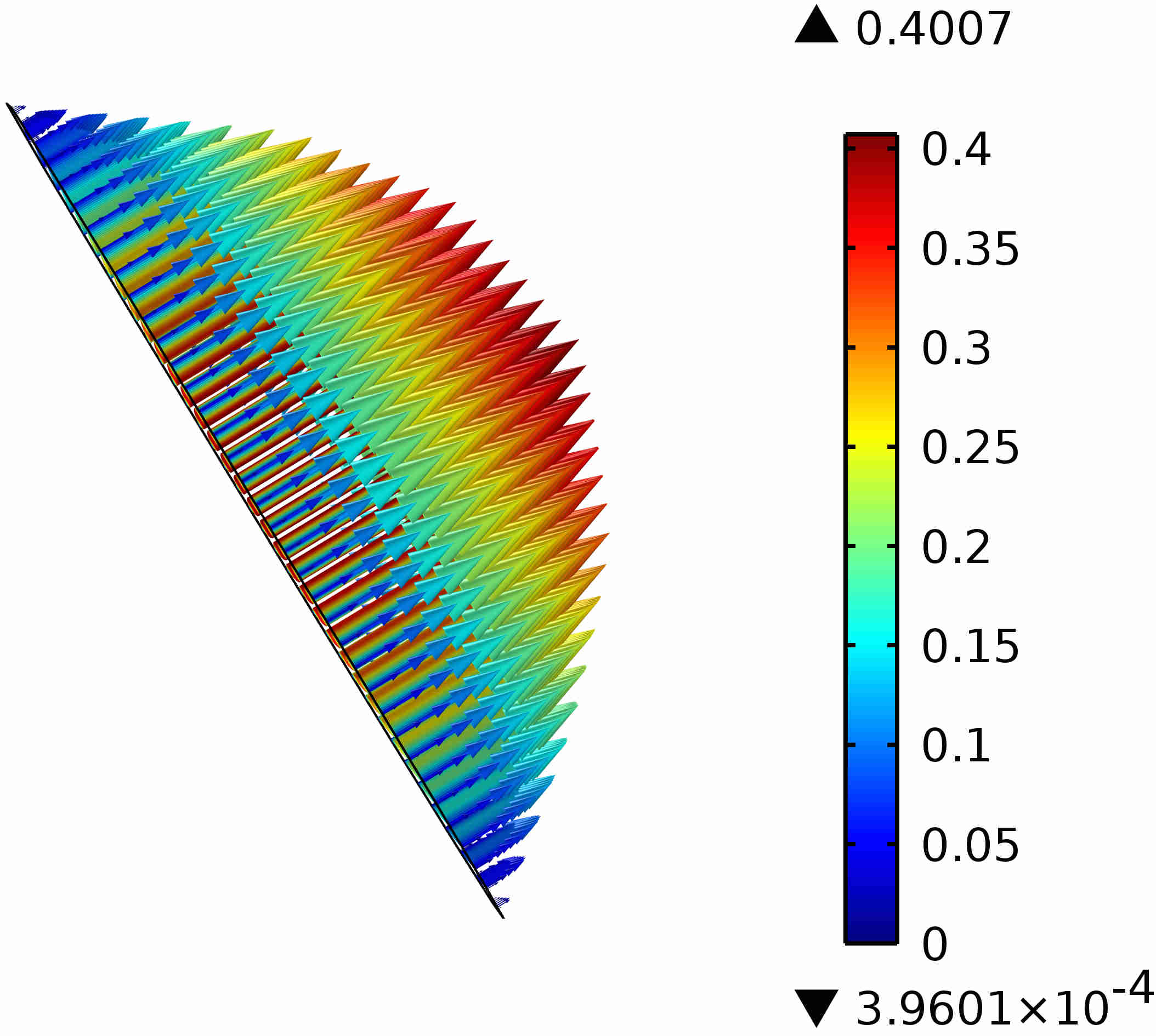}
\caption{Idealized   ${\bf u}_Q$.}
\end{subfigure}
\caption{Velocity vectors on $\Gamma_{in}$ {\it inlet 2}.}
\label{velcrime_in}
\end{figure}

To obtain the finite dimensional problem (\ref{nlcost})-(\ref{nlconst}) we used the same type of FEM and the same mesh; or, to be more precise, its part corresponding only to $\Omega$. As mentioned in beginning of section~\ref{realistic},  both the fact that we assume to know exactly the velocity at $\Omega_{part}$ and that we use the same mesh, put ourselves in the so called {\it inverse crime} scenario, but, in the frame of fluid control, such scenario is not necessarily trivial. Hence, before dropping this assumptions, we  ascertain that the DA approach can work - at least - in this case. The control problem was solved using Algorithm~\ref{SNOPT} with an optimality tolerance of $10^{-5}$. The simulation run for $1h15m$ on a Intel Xeon E5504 2.00 GHz using 4 cores. Looking to the first row of  Table~\ref{table3d1}, we can see the relative error of the controlled solution ${\bf u}$, with respect to the ground truth solution ${\bf u}_d$, evaluated at different parts of the domain. The relative error on $\Gamma_{in}$ gives us a measure of how the control vector differs from the ground truth solution at the artificial boundary. As we can confirm with the velocity profiles represented in Figure~\ref{velcrime_in} (b) and (a), respectively, they do not exactly match at $\Gamma_{in}$. Nevertheless, the relative errors $RE_{\Omega}$ and $RE_{\Omega_{part}}$ show a very good accuracy in the working  domain, and almost a perfect match at $\Omega_{part}$. 

\begin{figure}
\centering
 \begin{subfigure}[b]{0.25\textwidth}
\includegraphics[scale=0.2]{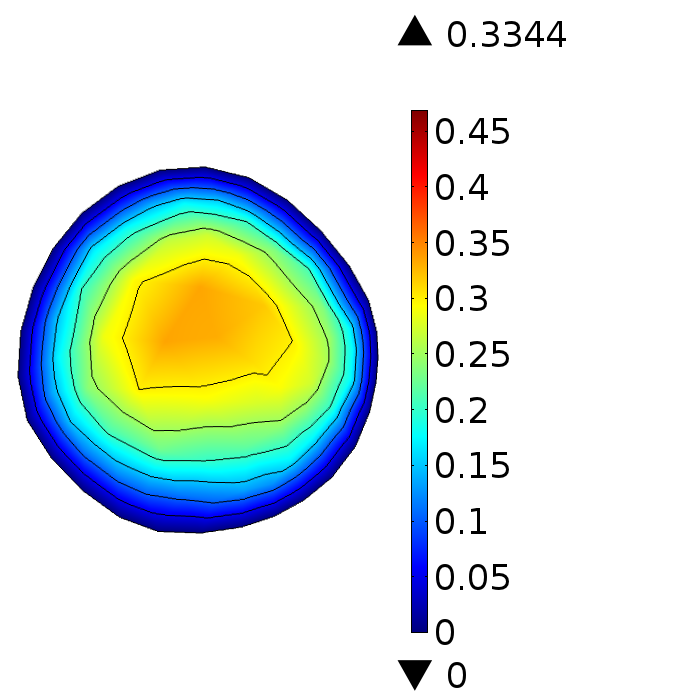}
\caption{S0: ${\bf u}_d$}   
\end{subfigure}
 \begin{subfigure}[b]{0.25\textwidth}
\includegraphics[scale=0.2]{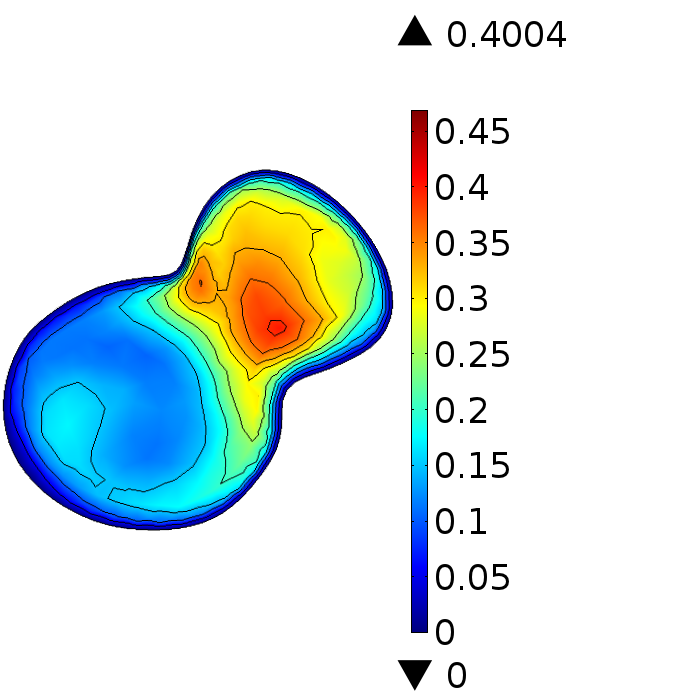}
\caption{S2:   ${\bf u}_d$}
\end{subfigure}
 \begin{subfigure}[b]{0.25\textwidth}
\includegraphics[scale=0.2]{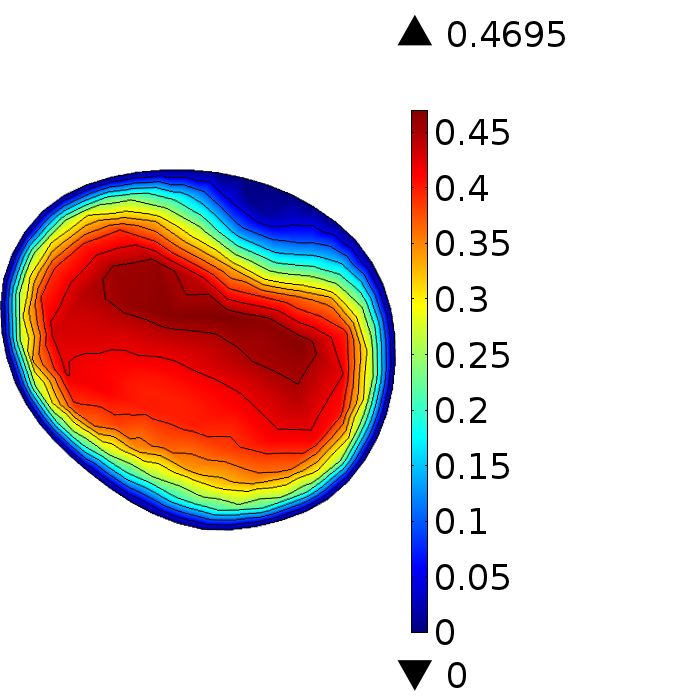}
\caption{S3:   ${\bf u}_d$}
\end{subfigure}
 \begin{subfigure}[b]{0.25\textwidth}
\includegraphics[scale=0.2]{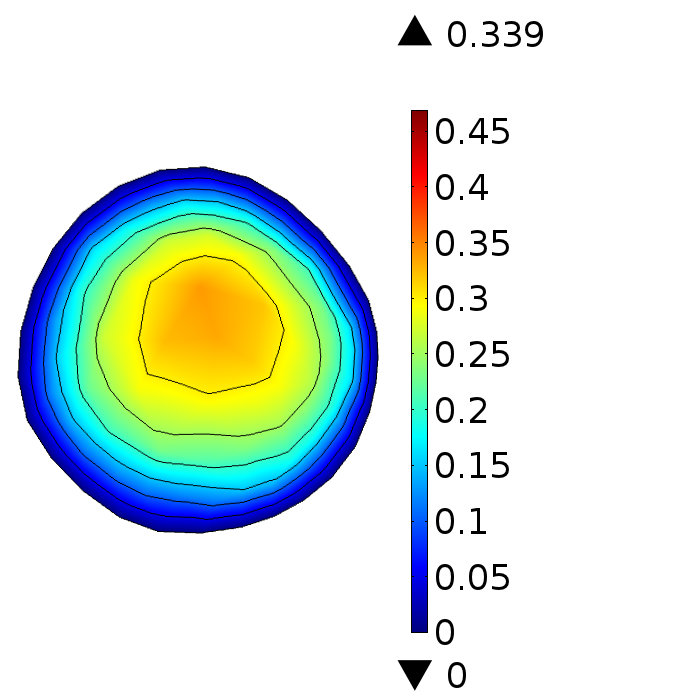}
\caption{S0: ${\bf u}$}
\end{subfigure}
 \begin{subfigure}[b]{0.25\textwidth}
\includegraphics[scale=0.2]{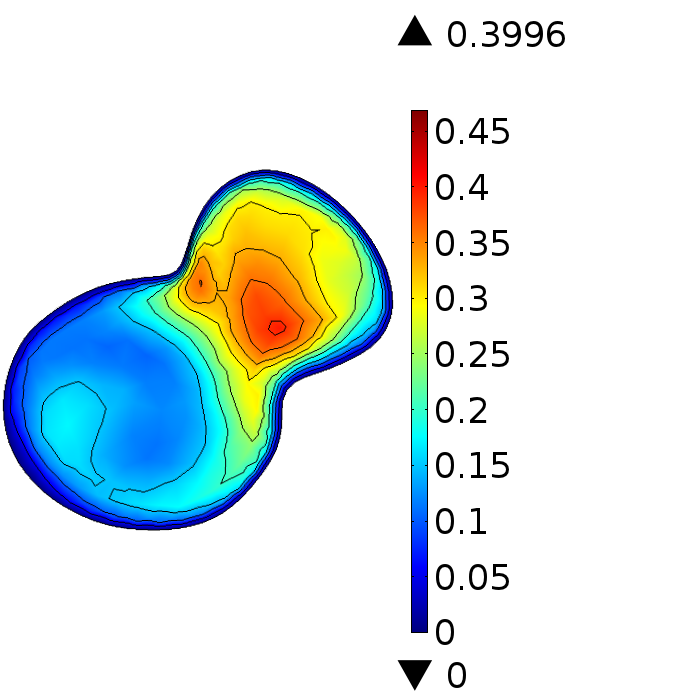}
\caption{S2:   ${\bf u}$}
\end{subfigure}
 \begin{subfigure}[b]{0.25\textwidth}
\includegraphics[scale=0.2]{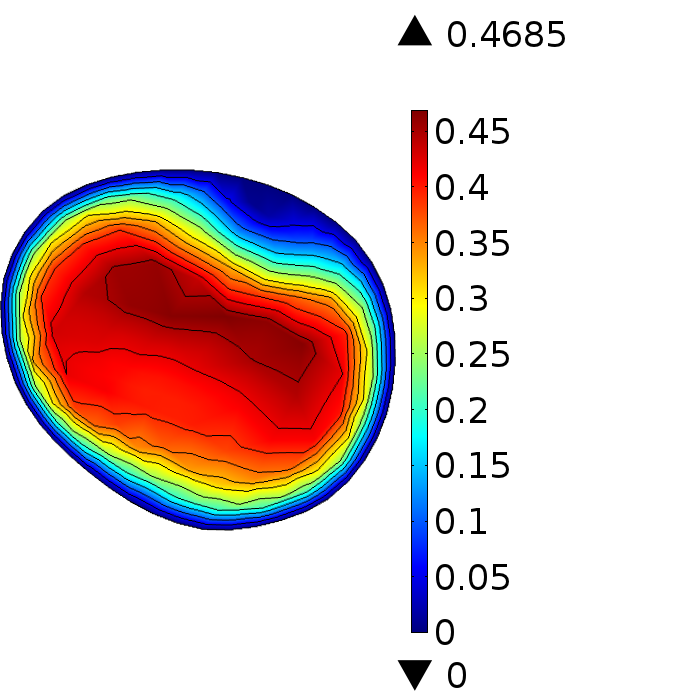}
\caption{S3:   ${\bf u}$}
\end{subfigure}
 \begin{subfigure}[b]{0.25\textwidth}
\includegraphics[scale=0.2]{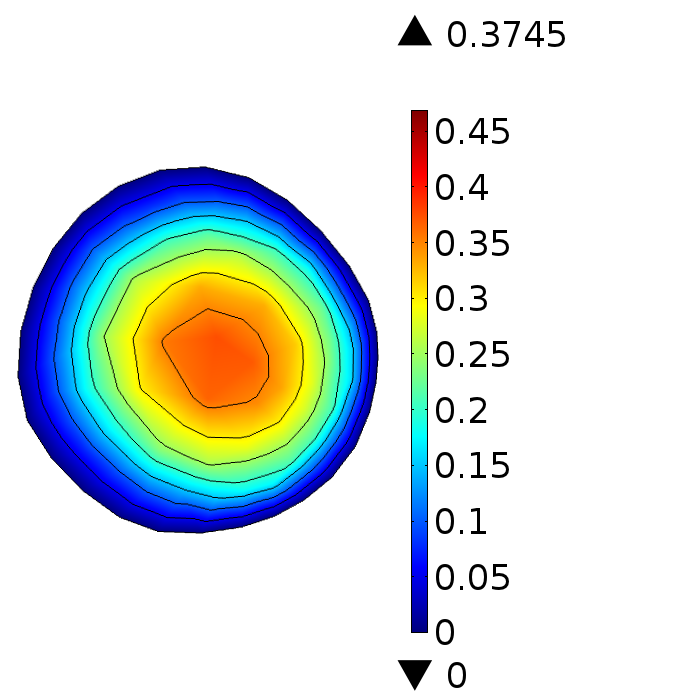}
\caption{S0: ${\bf u}_Q$}
\end{subfigure}
 \begin{subfigure}[b]{0.25\textwidth}
\includegraphics[scale=0.2]{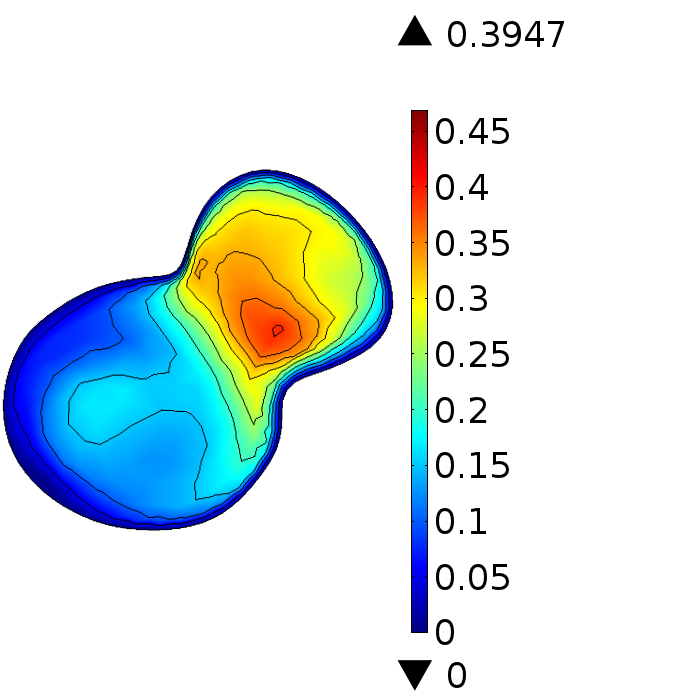}
\caption{S2:   ${\bf u}_Q$}
\end{subfigure}
 \begin{subfigure}[b]{0.25\textwidth}
\includegraphics[scale=0.2]{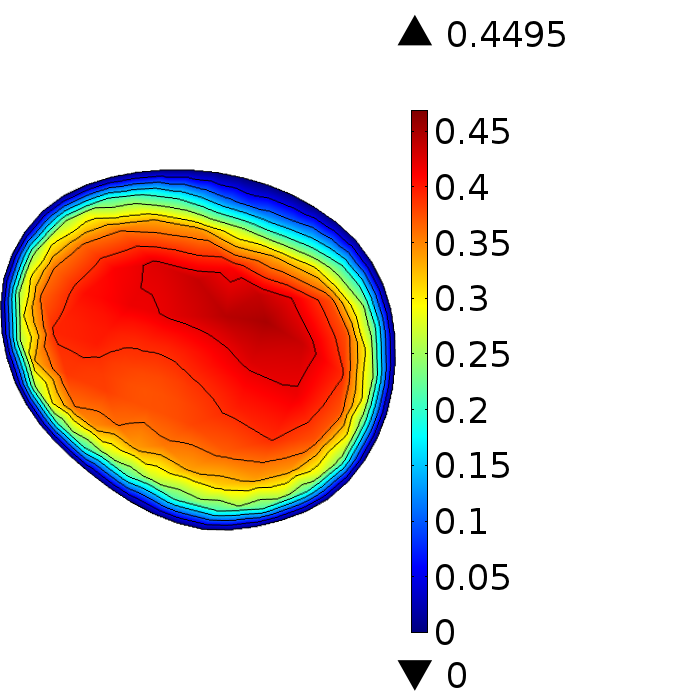}
\caption{S3:   ${\bf u}_Q$}         
\end{subfigure}
\caption{Velocity magnitude (m/s) for the ground truth ${\bf u}_d$, the controlled solution ${\bf u}$ and the idealized ${\bf u}_Q$ represented in several sections.}
\label{velcrime}
\end{figure}

\begin{figure}
\centering 
 \begin{subfigure}[b]{0.35\textwidth}
\includegraphics[scale=0.4]{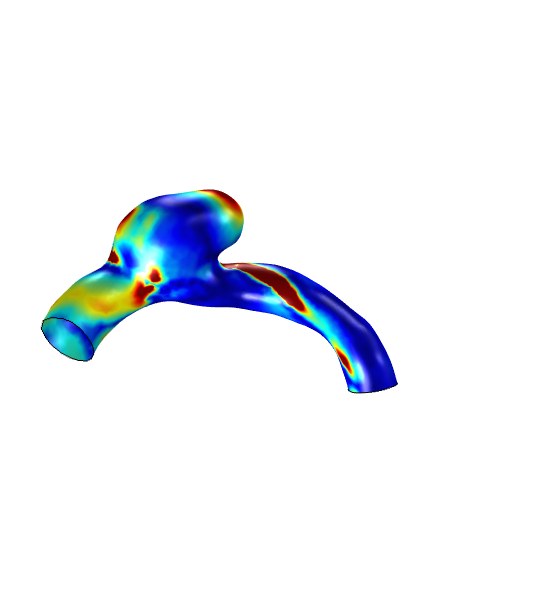}
\end{subfigure} 
 \begin{subfigure}[b]{0.45\textwidth}
\includegraphics[scale=0.4]{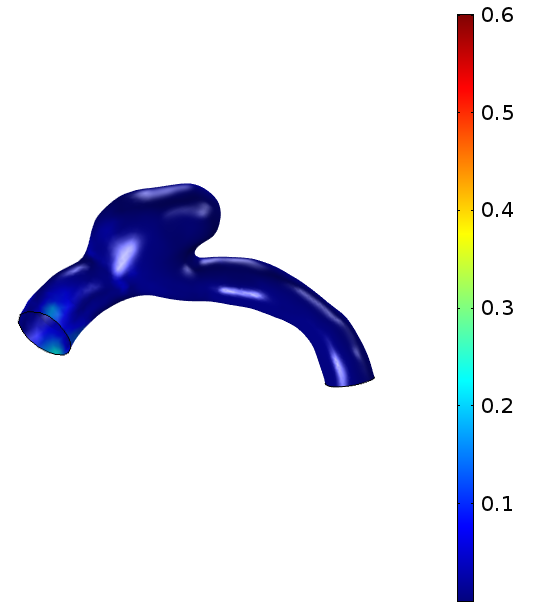}
\end{subfigure}
\caption{WSS ($N/m^2$): relative error with respect to ${\bf u}_d$. On the left, ${\bf u}_Q$.  On the right, ${\bf u}$. }
\label{WSSrelativeerror}
\end{figure}

To emphasize the gain achieved by the DA approach, we computed an alternative solution, ${\bf u}_Q$, based on the assumption that we can measure, instead of the velocity profile on $\Omega_{part}$, the exact flow rate. This rate was then used to define a laminar normal profile at $\Gamma_{in}$ (Figure~\ref{velcrime_in} (c)). The corresponding solution was done obtained similarly to ${\bf u}_d$. On the second row of Table~\ref{table3d1} we present the relative errors of this alternative solution with respect to the true solution ${\bf u}_d$. It can be found that the  DA approach - resulting in ${\bf u}$ - is prone for an error reduction, in the whole domain, from $22\%$ to less than $2\%$, when compared to the idealized solution ${\bf u}_Q$. In Figure~\ref{velcrime} a representation of  the   velocity magnitudes for the true solution ${\bf u}_d$, the controlled solution ${\bf u}$ and the idealized ${\bf u}_Q$ is depicted along three sections. The similarities between ${\bf u}$ and ${\bf u}_d$ are clear. In Figure~\ref{WSSrelativeerror} we can see the representation of the relative error, for the WSS magnitude, of both ${\bf u}$ (right) and ${\bf u}_Q$ (left) with respect to the ground truth solution ${\bf u}_d$. While for the solution ${\bf u}_Q$, based on the idealized laminar profile, the relative error is frequently above $40\%$ and sometimes above $60\%$, the relative error associated to the controlled solution only reaches $10\%$ close to the inlet, where the mesh was chosen deliberately coarse. This indicates an important potential gain of the DA approach in reducing the error associated to WSS {\it in silico} measurements.

\begin{table}[ht!]
\caption{Relative errors and final value for cost functional (J) for (P1) in the realistic domain. }
\centering

\begin{tabular}{ccccc}
\toprule
Solution  & $RE_{\Omega}$ &$RE_{\Gamma_{in}}$ & $RE_{\Omega_{part}}$ & Cost\\
\midrule
${\bf u}$   & $0.019522$ & $0.176382$ & $0.002823$ & $0.001037$  \\ 
${\bf u}_Q$     & $0.225058$ & $0.40286$ & $0.25353$ & $0.04693$\\
\bottomrule
		  \end{tabular}
		\label{table3d1}
\end{table}

\subsubsection{Generating Data with a finer mesh.}
 The results that we present next correspond to the case where the ground truth solution ${\bf u}_d$ is not generated on the same mesh used for the DA approach (the control problem). Specifically, the later is conducted using the same mesh, as in the previous example, while the former is obtained using 412K degrees of freedom. In this way, we dropped one of the two assumptions placing us in the {\it inverse crime} scenario. 
 
 We remark that the relative error of the ground truth solution from the previous example (213K), with respect to the new ${\bf u}_d$ computed in this way ($412K$), is approximately $14\%$, even if the same exact boundary conditions are used. As a consequence, to obtain a better match inside the domain, one must expect that the boundary conditions will not exactly coincide. Once more we assume to know the velocity at $\Omega_{part}=S1\cup S2\cup S3\cup S4$ but, this time,  in the form of a linear interpolation of the ground truth data. In the first row of Table \ref{table3d2} we can see the relative errors and final cost functional resulting from  the DA approach in this scenario which we refer to as (P3). We see that, even if the control doesn't match the data at the artificial boundary $\Gamma_{in}$, but allows  a $50\%$ reduction in $RE_{Omega}$ when compared with the error measured for ${\bf u}_Q$.  Also, the relative error on $\Omega_{part}$ is around $6\%$, one quarter of the error associated to ${\bf u}_Q$. In fact, since the solutions obtained for the two meshes do not coincide (even if imposing the same boundary conditions), minimizing (\ref{costfunctional}) forces the control to differ from the data on the inlet, in order to accomplish an improved matching on $\Omega_{part}$.  Looking at the  first two rows of Figure~\ref{velmag4},  we can see that the controlled solution ${\bf u}$ is accurate on the sections that belong to $\Omega_{part}$ (i.e., $S2$ and $S3$), while it distinguishes from ${\bf u}_d$ in $S0$, the section close to the artificial inlet. This finding is consistent with the previous comments.

{To check the robustness of the approach with respect to increasing Reynolds numbers, we repeated the simulations the scenario just described, using different Reynolds numbers. In Table \ref{table  ReynoldsIII} the results are shown. The relative errors on the observations remain of the same order. We remark that these results can not be directly compared with the results from the previous 2D example, as the SNOPT library automatically adjusts for the different characteristics of the underlying optimization problem.}

\begin{table}[ht!]
\caption{Relative errors, number of SNOPT iterations for different Reynolds numbers. }
\centering

\begin{tabular}{cccccc}
\toprule
Reynolds  & $RE_{\Omega}$ &$RE_{\Gamma_{in}}$ & $RE_{\Omega_{part}}$ & &Iterations\\
\midrule
$367$    & $0.106451$ & $0.500379$ & $0.067542$ & $$&$372$\\
$500$   & $0.104979$ & $0.492143$ & $0.067777$ & $$&$364$\\
$1000$   & $0.104229$ & $0.481692$ & $0.068134$ & $$&$363$\\
$2000$   & $0.105774$ & $0.495097$ & $0.067668$ & $$&$370$\\
\bottomrule
\end{tabular}
\label{table ReynoldsIII}
\end{table}

\subsubsection{Adding noise to the observations.}\label{secnocrime}

 In real life examples,  data measurements frequently include  a certain error in the form of noise, due to the lack of accuracy of observation devices. As it has been emphasized in \cite{Delia3} and \cite{GTS}, the DA variational approach has  an important role in noise reduction. To mimic this scenario, we now consider the case  where the data available at $\Omega_{part}$ is perturbed with noise.  To represent such noise, we randomly generated a sample from a normal distribution with zero mean and standard deviation given by $\bar{\sigma}=0.2\frac{U_0}{3}$, where $U_0=0.419507\,m/s$ is the maximum velocity of ${\bf u}_d$ at the inlet boundary. We then computed ${\bf u}$ as above. As we can see on the second row of Table \ref{table3d2}, the presence of noise is overcome, and an accurate solution, similar to the one obtained in the case without noise, is obtained at $\Omega_{part}$. This can be confirmed from the pictures represented on the first to third rows  of Figure~\ref{velmag4}. These results are in agreement with the conclusions of \cite{Delia3} and \cite{GTS}. Additionally, it can be also be seen (5th row) how the idealized solution ${\bf u}_Q$ differs more from the data, when compared to the controlled solutions.

{With regards with the choice of the parameters $\beta_1$ and $\beta_2$, we proceeded similarly to the previous example but, instead of a reference relative error,  we used a reference noise $L^2$ norm corresponding to a magnitude of order $\delta\approx 10^{-4}$. We tested a sample of parameters. Among several possibilities verifying $\|{\bf u}_{\beta_1,\beta_2}-{\bf u}_d\|_{\mathbf{L}^2(\Omega_{part})}\approx 10^{-4}$, both choices $(\beta_1,\beta_2)_1=(0.7144\times 10^{5},10^{-3})$ and $(\beta_1,\beta_2)_2=(10^{5},10^{-3})$ resulted in similar relative errors, of $RE_{\Omega_{part}}=0.0646$ and $RE_{\Omega_{part}}=0.0675$, respectively. The first option was considered because it allowed a normalization of the first term in the cost function. However, it was the second choice that allowed to  computational solve all the test cases here shown, using exactly the same solver configurations. Therefore, for the realistic example, we fixed $(\beta_1,\beta_2)=(10^{5},10^{-3})$}.

\subsubsection{Reducing the observed set.}
\begin{figure}
\centering
 \begin{subfigure}[b]{0.25\textwidth}
\includegraphics[scale=0.2]{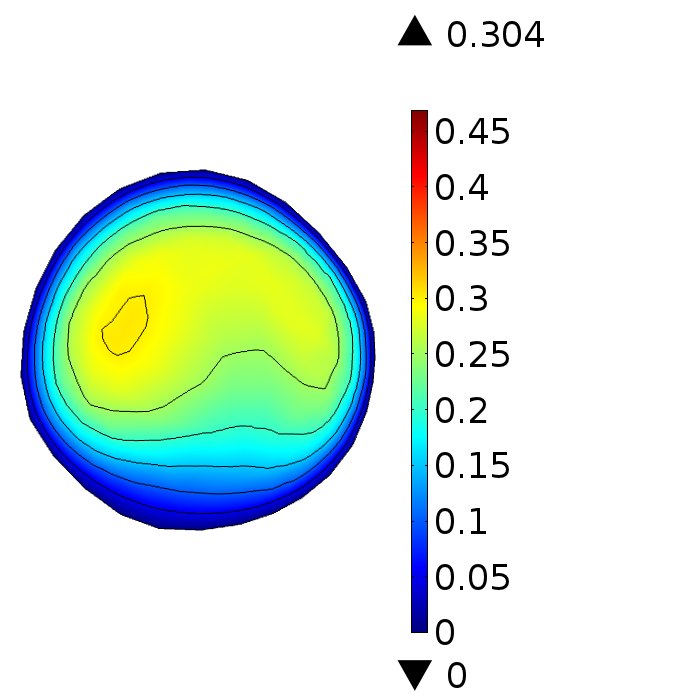}
\caption{S0:  imported ${\bf u}_d$}
\end{subfigure}
 \begin{subfigure}[b]{0.25\textwidth}
\includegraphics[scale=0.2]{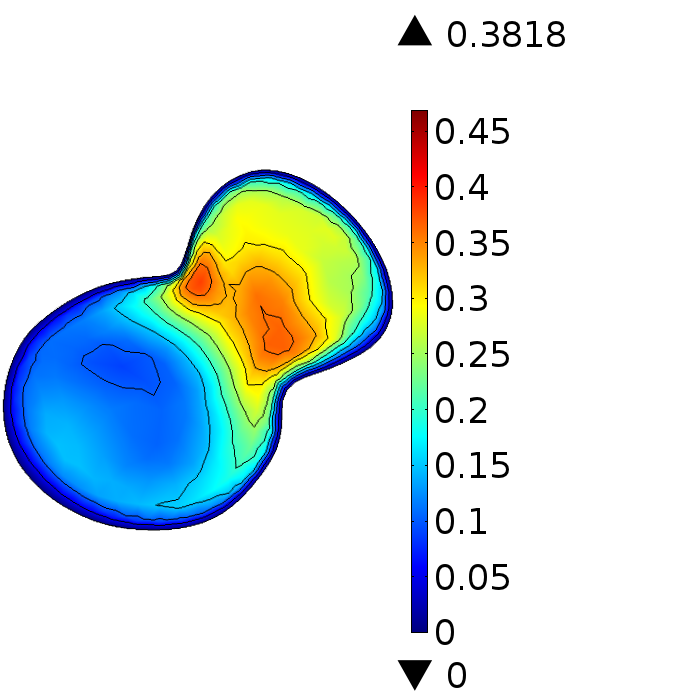}
\caption{S2:  imported ${\bf u}_d$}
\end{subfigure}
 \begin{subfigure}[b]{0.25\textwidth}
\includegraphics[scale=0.2]{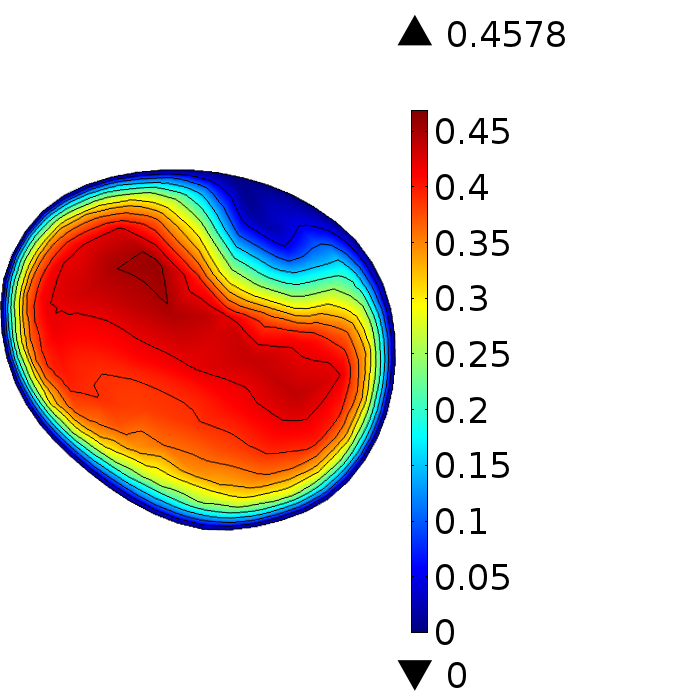}
\caption{S3:  imported ${\bf u}_d$}
\end{subfigure}
 \begin{subfigure}[b]{0.25\textwidth}
\includegraphics[scale=0.2]{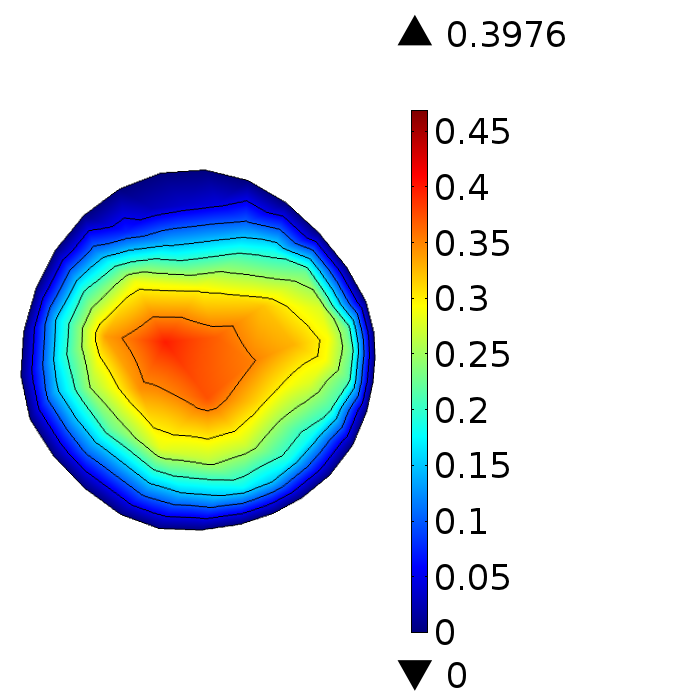}
\caption{S0:  P3 }
\end{subfigure}
 \begin{subfigure}[b]{0.25\textwidth}
\includegraphics[scale=0.2]{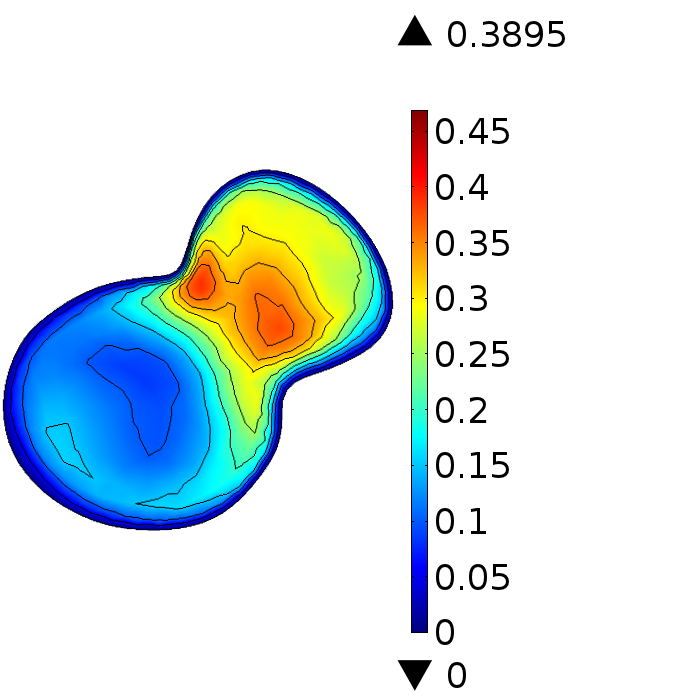}
\caption{S2:   P3 }
\end{subfigure}
 \begin{subfigure}[b]{0.25\textwidth}
\includegraphics[scale=0.2]{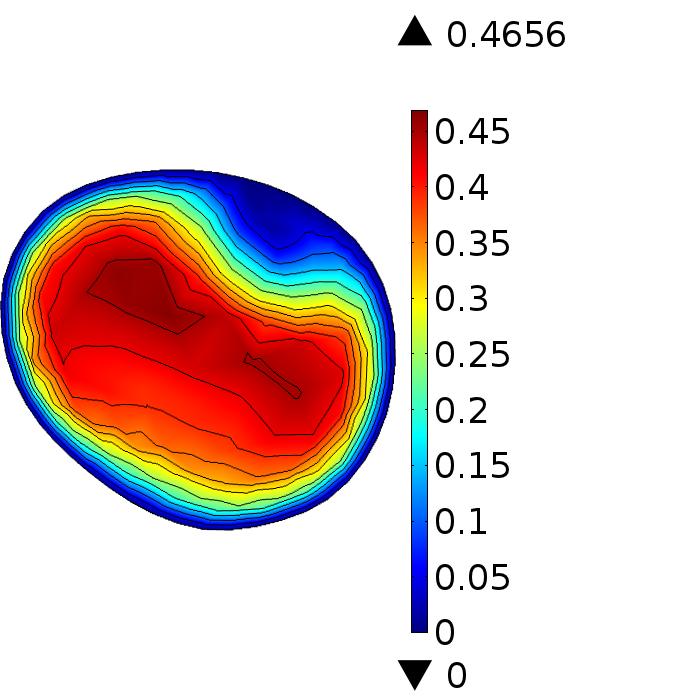}
\caption{S3:    P3 }
\end{subfigure}
 \begin{subfigure}[b]{0.25\textwidth}
\includegraphics[scale=0.2]{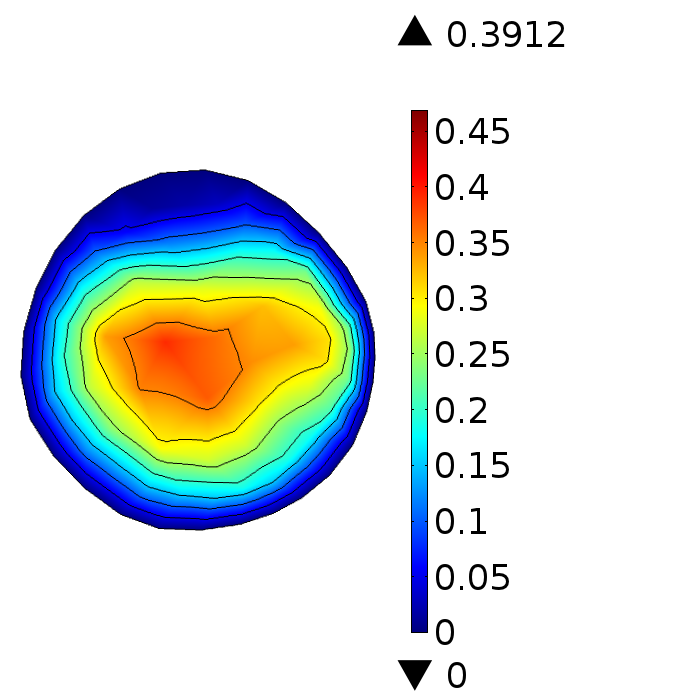}
\caption{S0: Noisy P3 }
\end{subfigure}
 \begin{subfigure}[b]{0.25\textwidth}
\includegraphics[scale=0.2]{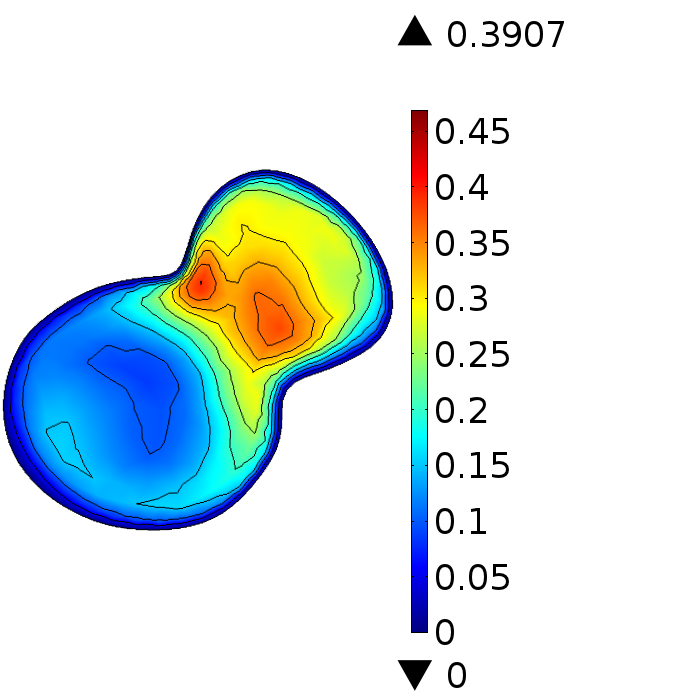}
\caption{S2: Noisy P3 }
\end{subfigure}
 \begin{subfigure}[b]{0.25\textwidth}
\includegraphics[scale=0.2]{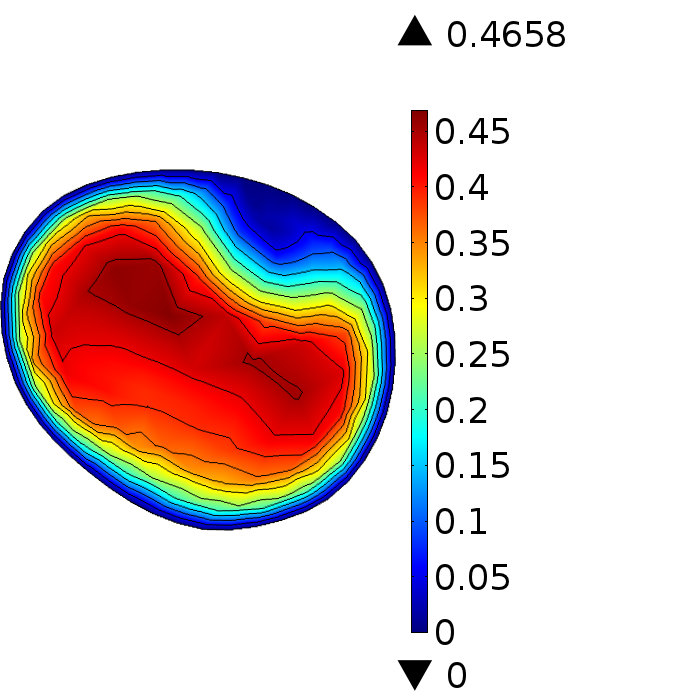}
\caption{S3: Noisy P3 }
\end{subfigure}
 \begin{subfigure}[b]{0.25\textwidth}
\includegraphics[scale=0.2]{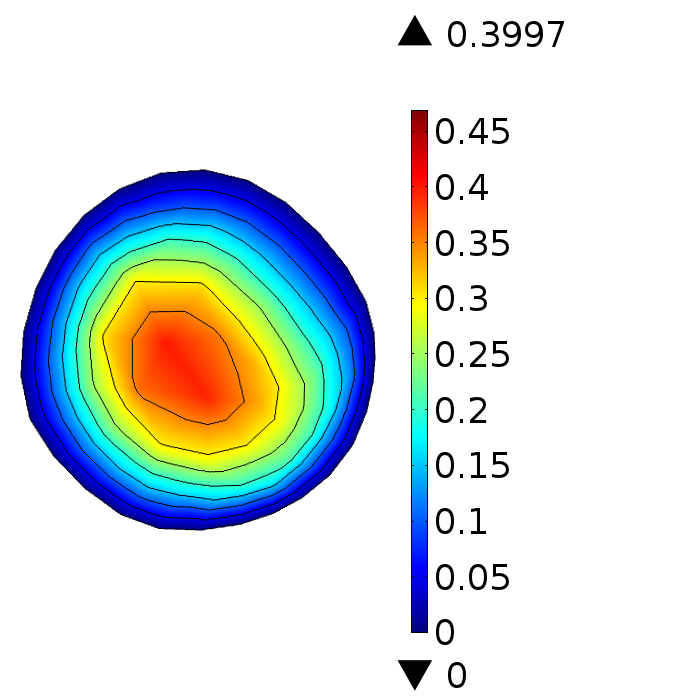}
\caption{S0:  P4}
\end{subfigure}
 \begin{subfigure}[b]{0.25\textwidth}
\includegraphics[scale=0.2]{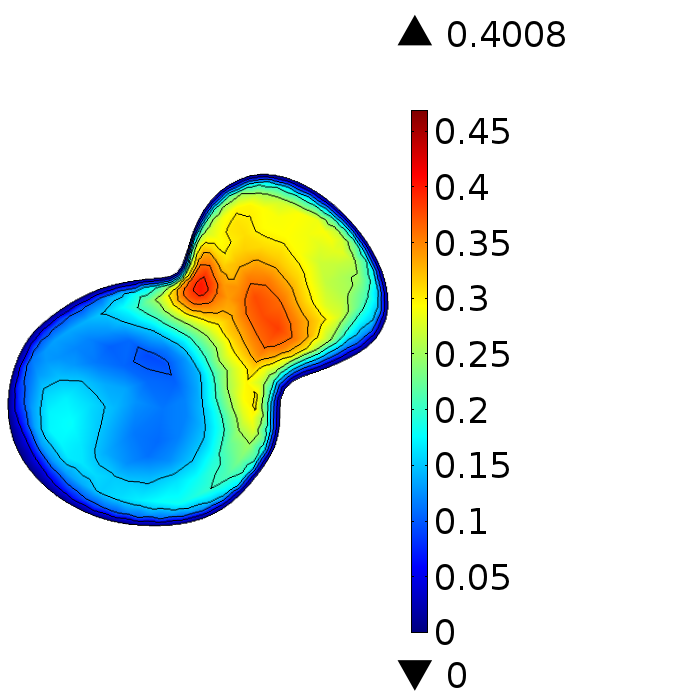}
\caption{S2:  P4}
\end{subfigure}
 \begin{subfigure}[b]{0.25\textwidth}
 \includegraphics[scale=0.2]{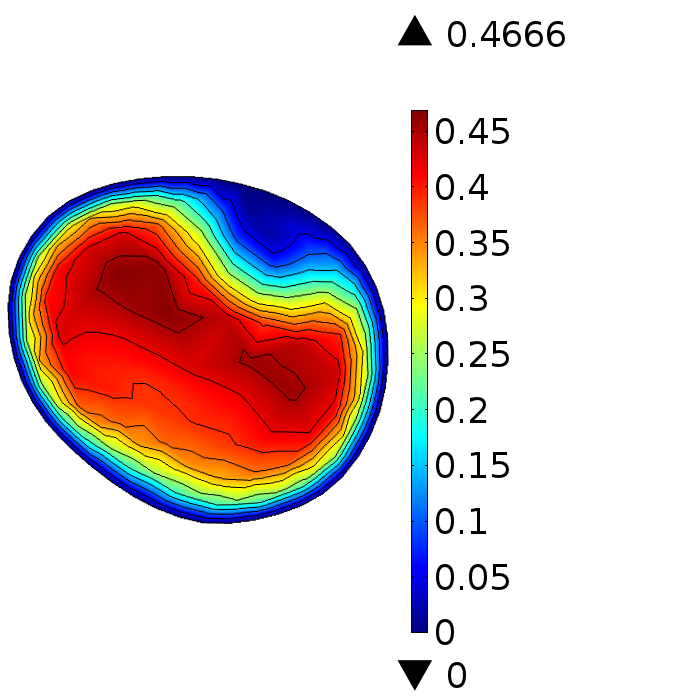}
\caption{S3:  P4}
\end{subfigure}
\begin{subfigure}[b]{0.25\textwidth}
\includegraphics[scale=0.2]{2d_vel_sec0_uQ.png}
\caption{S0: ${\bf u}_Q$}
\end{subfigure}
 \begin{subfigure}[b]{0.25\textwidth}
\includegraphics[scale=0.2]{2d_vel_sec2_uQ.png}
\caption{S2:   ${\bf u}_Q$}
\end{subfigure}
 \begin{subfigure}[b]{0.25\textwidth}
\includegraphics[scale=0.2]{2d_vel_sec3_uQ.png}
\caption{S3:   ${\bf u}_Q$}         
\end{subfigure}
\caption{ No crime scenario. Velocity magnitude (m/s) in several sections.}
\label{velmag4}
\end{figure}

Finally, in order to understand the role of the sections chosen to integrate $\Omega_{part}$, we consider the case where $\Omega_{part}=S3$, that is, we assume to have measurements only at the section located immediately downstream the aneurysm.  In fact, as mentioned in Section~\ref{intro}, velocity measurements inside the lumen are not trivially obtained, so it is actually more realistic to assume that we have fewer observations available. We refer to this case as (P4). On the third row of Table~\ref{table3d2} we can verify that the overall error increases from $10.6\%$ (4 sections) to $15.5\%$ (1 section), while the error in $\Omega_{part}=S3$  slightly decreases. The final value of the cost functional is smaller than for (P3), which is natural, since it consists of just a part of the cost in (P3). On the fourth row of Figure~\ref{velmag4}, we can realize that, whilst the controlled solution remains accurate at $S3$, it becomes more distinct from ${\bf u}_d$ when looking further upstream, at $S2$ and  $S0$. These findings agree with  the intuitive idea that extended measurements along the domain improve the overall accuracy of the controlled solution. A mathematical result to ensure these principle should rely on the observability and controllability concepts and should be treated more carefully, in a future work.

\begin{table}[ht!]
  \caption{Relative errors of the controlled solutions obtained by solving (P3) with and without noise, (P4) and for comparison purposes, of ${\bf u}_Q$.} 
\centering

\begin{tabular}{ccccc}
\toprule
Problem  & $RE_{\Omega}$ &$RE_{\Gamma_{in}}$ & $RE_{\Omega_{part}}$ & Cost\\
\midrule
(P3)  & $0.106451$ & $0.500379$ & $0.067542$ & $0.004842$\\

(P3) + noise & $0.110565$ & $0.526449$ & $0.067358$ & $0.037602$\\

(P4)   & $0.15527$ & $0.44471$ & $0.05703$ & $0.00149$\\

${\bf u}_Q$  & $0.225058$ & $0.40286$ & $0.25353$ & $0.04693$\\  
\bottomrule
		  \end{tabular}
		\label{table3d2}
\end{table}

\section{Conclusions}\label{conclusions}
In this work we have suggested a  velocity control approach, as a Data Assimilation (DA) technique in the frame of blood flow simulations. We have shown that the nonlinear control problem, inherent to this approach, is mathematically sound at the continuous level. A discretize then optimize procedure, followed by the application of a large scale sequential quadratic programming implementation, resumes the methodology.

By applying the approach to a suitable idealized example, we have identified the potential advantages with respect to a pressure control strategy. This is not a general conclusion, but it is valid in some cases, when the ground truth pressure profiles are not parallel to domain cross sections. The next step consisted in analyzing  a realistic situation where the computational domain was obtained from medical images of a brain aneurysm. Due to the lack of real velocity measurements, {\it in silico} profiles were generated to supply the required data. The results are promising even when high extra noise is added to the data: the error relative to the ground true solution is significantly reduced, when compared to a solution obtained from a laminar idealized profile.

At the present stage, the major drawback concerns the computational cost. In fact, a reliable WSS estimate requires 5 to 10 times more degrees of freedom. Also, as it was mentioned before, a definitive approach should allow the coupling with extra models describing fluid structure interaction and certain pathologies, such as clot formation or plaque growth. Finally, a straightforward application of this methodology to a time dependent simulation would have this computational cost associated to each time step iteration. Therefore, even if for the first time, DA techniques were validated as a means to recover fully general velocity profiles, further improvements must consider the inclusion of order reduction techniques, such as the ones suggested in \cite{BV}. Nevertheless, the authors believe that the results here are the necessary sound basis for such future improvements.


\pagebreak


\end{document}